\input  amstex
\input amsppt.sty
\magnification1200
\vsize=23.5truecm
\hsize=16.5truecm
\NoBlackBoxes

\def\supp{\operatorname{supp}}

\def\crp{\overline{\Bbb R}_+}
\def\crm{\overline{\Bbb R}_-}
\def\crpm{\overline{\Bbb R}_\pm}

\def\rnp{{\Bbb R}^n_+}
\def\rnm{\Bbb R^n_-}
\def\rnpm{\Bbb R^n_\pm}
\def\crnp{\overline{\Bbb R}^n_+}
\def\crnm{\overline{\Bbb R}^n_-}

\def\comega{\overline\Omega }
\def\ang#1{\langle {#1} \rangle}
\def\simto{\overset\sim\to\rightarrow}

\def\OP{\operatorname{OP}}
\def\rp{ \Bbb R_+}

\def\rpm{ \Bbb R_\pm}

\def\R{\Bbb R}
\def\C{\Bbb C}

\def\ol{\overline}
\def\SD{\Cal S}
\def\F{\Cal F}
\def\E{\Cal E}
\def\D{\Cal D}
\def\oc{\overset\sssize\circ\to}

\document

\medskip

\topmatter
\title
Fractional Laplacians on domains,\\ 
a development of H\"o{}rmander's theory of\\ mu-transmission
pseudodifferential operators
\endtitle
\author Gerd Grubb \endauthor
\affil
{Department of Mathematical Sciences, Copenhagen University,
Universitetsparken 5, DK-2100 Copenhagen, Denmark.}
E-mail {\tt grubb\@math.ku.dk}
\endaffil
\rightheadtext{Fractional Laplacians}
\abstract
Let $P$ be a classical pseudodifferential operator of order $m\in{\Bbb
C}$ on an $n$-dimensional $C^\infty $ manifold $\Omega _1$. For the
truncation $P_\Omega $ to a smooth subset $\Omega $  there is a well-known theory of boundary value problems when $P_\Omega $ has the
transmission property (preserves $C^\infty (\comega)$) and is of
integer order; the calculus of Boutet de Monvel. Many interesting
operators, such as for example complex powers of the
Laplacian $(-\Delta )^\mu $ with $\mu \notin{\Bbb Z}$, are not covered. They
have instead the $\mu $-transmission property defined in H\"ormander's books, mapping $x_n^\mu C^\infty (\comega)$ into $C^\infty (\comega)$.
In an unpublished lecture note from 1965,
H\"ormander described an $L_2$-solvability theory
for $\mu $-transmission  operators, departing from Vishik and Eskin's results.  
We here
develop the theory in $L_p$ Sobolev spaces ($1<p<\infty $) in a modern setting.
It leads to not only Fredholm solvability statements but also regularity results
in full scales of Sobolev spaces ($s\to\infty $). The solution spaces
have a singularity at the boundary that we describe in detail. We moreover obtain results in H\"older spaces, which radically improve
recent regularity results for fractional Laplacians. 
 
\endabstract
\dedicatory To the memory of Lars H\"ormander 1931--2012
\enddedicatory
\endtopmatter

\subhead  Introduction \endsubhead
Pseudodifferential operators ($\psi $do's) of integer order with the transmission property
(preserving $C^\infty $ up to the boundary in a domain)
and their boundary problems have been studied since the basic
theory was developed by Boutet de Monvel in \cite{B71}. The theory includes differential operators and
the parametrices of elliptic such ones, and also operators whose
symbols are rational functions of $\xi $. 

This was preceded by works of Vishik and Eskin (\cite{VE65},
\cite{VE67} etc., included for the
major part in Eskin's book \cite{E81}), which treated operators of a more
general type, having a factorization of the principal symbol at the
boundary of a smooth open set $\Omega $, in two
factors extending analytically to $\{\operatorname{Im}\xi _n>0\}$
resp.\ $\{\operatorname{Im}\xi _n<0\}$ as functions of the conormal
variable $\xi _n$, with each their degree of homogeneity $m-\kappa
(x')$ resp. $\kappa (x')$, $x'\in\partial\Omega $.
When $\Omega $ is compact, such operators will under mild restrictions
on the factorization index $\kappa (x')$ define Fredholm operators on
Sobolev spaces with exponent $s$ in a certain open interval
$\,]s_-,s_+[\,$ of length $\le 1$. For larger $s$ one has to
 add suitable boundary
conditions, and for smaller $s$ potential terms, in order to get
Fredholmness.
The results have been extended to $L_p$-based Sobolev
spaces by Shargorodsky \cite{S94} and 
Chkadua and Duduchava  \cite{CD01}.

In an unpublished (photocopy distributed) lecture note
 at Princeton 1965  \cite{H65},
H\"o{}r\-man\-der introduced, with Vishik and Eskin's work as a starting point,
 a generalized transmisssion condition of type $\mu \in{\Bbb C}$ (where
the condition in \cite{B71} is the case $\mu =0$), reflecting the
properties of the general operators studied by Vishik and Eskin in the
 case $\kappa (x')=\mu _0$ constant. Here he showed
not only the Fredholm property in Sobolev spaces for $s$ in an
interval, but he moreover determined the $L_2$ Sobolev regularity of solutions with
data given for all larger $s$, or given in $C^\infty (\comega
)$, finding the domain spaces for Fredholm solvability and describing
 the associated boundary conditions.
 
The transmission condition of type $\mu $ was briefly characterized in
\cite{H85}, Sect.\ 18.2. An application to 
propagation of singularities was given by Hirschowitz and Piriou \cite{HP79}.

Fractional powers of the Laplacian $(-\Delta )^a$ are of type
$\mu =a$; they have recently received increased attention both in probability theory, cf.\ e.g.\
Bogdan, Grzywny and Ryznar \cite{BGR10}, 
Ros-Oton and Serra
\cite{RS14}, in differential geometry, cf.\ e.g.\ Gonzalez, Mazzeo and
Sire \cite{GMS12}, and in Schr\"odinger theory, cf.\
e.g.\ Frank and Geisinger \cite{FG14}, and the references in these
papers. Only a little seems to be known about the regularity of
solutions on domains. Inspired by this, we have in the present paper worked out an
extension of H\"ormander's theory to
$L_p$-Sobolev spaces, $1<p<\infty $, with additional results, moreover leading to solvability results in
H\"older spaces. Applications include
fractional powers of strongly elliptic differential operators.

In this process, the presentation could benefit from the theories
developed since 1965, namely the
theory of boundary value problems of type 0, as introduced by Boutet
de Monvel for integer-order cases in \cite{B71}, and further developed
 by the present author, e.g.\ in \cite{G96}. The work \cite{G90} is
 particularly useful, extending the Boutet de Monvel calculus to the
$L_p$-setting and introducing refined order-reduction techniques. A
joint work with H\"ormander \cite{GH90} treated operators of type 0 and arbitrary
real order $m$ (including $S^m_{\varrho ,\delta }$ symbols).

Here are some of the main results.
We consider
a smooth subset $\Omega $ of an
$n$-dimensional Riemannian $C^\infty $ manifold $\Omega _1$, and
denote by $d(x)$ a $C^\infty (\comega)$-function equal to
$\operatorname{dist}(x,\partial\Omega )$ near $\partial\Omega $ and
positive on $\Omega $. Restriction to $\Omega $ is denoted $r_\Omega $ (or $r^+$), extension
by zero on $\Omega _1\setminus \Omega $ is denoted $e_\Omega $ (or
$e^+$). 
For $\mu
\in {\Bbb C}$ with $\operatorname{Re}\mu >-1$,  $\Cal E_\mu (\comega)$
denotes the space of functions $u$ such that
$u=e_\Omega d(x)^\mu v$ with $v\in C^\infty (\comega)$. The definition is 
generalized in a distribution sense to lower values of $\mu $.
On $\Omega _1$ we consider  a classical $\psi $do
$P$ of order $m\in{\Bbb C} $,
with symbol in local coordinates $p(x,\xi )\sim \sum_{j\in{\Bbb N}_0}p_j(x,\xi )$ where
$p_j(x,t\xi )=t^{m -j}p_j(x,\xi )$. The $\mu $-transmission property
was described in \cite{H85}, Th.\ 18.2.18:

\proclaim{Proposition 1} A necessary and sufficient condition in order
that $r_\Omega Pu\in C^\infty (\comega)$ for all $u\in \E_\mu
(\comega)$ is that $P$ satisfies the $\mu $-transmission condition (in
short: is of type $\mu $), namely  that $$
\partial_x^\beta \partial_\xi ^\alpha {p_j}(x,-N)=e^{\pi i(m-2\mu -j-|\alpha | )
}\partial_x^\beta \partial_\xi ^\alpha{p_j}(x,N),\; x\in\partial\Omega
,\tag1
$$
for all $j,\alpha ,\beta $, where 
$N$ denotes the
interior normal to $\partial\Omega $ at $x$.
\endproclaim

In the following
theorems we take $\comega $ compact. 

Define the special spaces $H_p^{\mu (s)}(\crnp)$ (H\"ormander's $\mu $-spaces), for $s>\operatorname{Re}\mu -1/p'$:
$$
 H_p^{\mu (s)}(\crnp)=\{u\in
\dot H_p^{\operatorname{Re}\mu -1/p'+0}(\crnp)\mid r^+\OP((\ang{\xi '}+i\xi _n)^\mu )u\in
\ol H_p^{s-\operatorname{Re}\mu }(\rnp)\}.\tag2
$$
(The notation used for $L_p$ Sobolev spaces is listed below in Section 1.)
The definition extends to define $H_p^{\mu (s)}(\comega)$ by use of
local coordinates. This is the solution space for $Pu=f$ on $\Omega $:

\proclaim{Theorem 2} Assume that $P$ is elliptic
 of order $m\in{\Bbb C}$ and type $\mu _0\in{\Bbb C}$ {\rm (mod 1)}, and has factorization index $\mu _0$, and let $s>\operatorname{Re}\mu _0-1/p'$.
When  $u\in\dot H_p^{\operatorname{Re}\mu _0-1/p'+0}(\comega)$, then $r_\Omega Pu\in \ol
H_p^{s-\operatorname{Re}m }(\Omega )$ implies $u\in H_p^{\mu _0(s)}(\comega)$. 
The mapping $$
r_\Omega P\colon H_p^{\mu _0(s)}(\comega)\to  \ol
H_p^{s-\operatorname{Re}m }(\Omega )\tag3$$ 
is Fredholm.
Moreover, $r_{\Omega }Pu\in C^\infty (\comega)$ implies   $u\in
\E_{\mu _0} (\comega)$, and
the mapping $r_\Omega P$
from $\E_{\mu _0} (\comega)$ to $C^\infty (\comega)$ is Fredholm.
\endproclaim

The spaces $H_p^{\mu (s)}(\comega )$ allow a definition of 
{\it boundary values} $\gamma _{\mu ,j}u$, that generalize the
mapping $u\mapsto \partial_{x_n}^j(x_n^{-\mu
}u)|_{x_n=0}$, defined for $u\in\E_{\mu }(\crnp)$
when $\operatorname{Re}\mu >-1$. 

\proclaim{Theorem 3} When $P$ and $s$ are as in Theorem {\rm 2},
and $\mu =\mu _0-M$ for a positive integer $M$, then
the following operator is Fredholm:
$$
\{r_\Omega P,\gamma _{\mu ,0},\dots,\gamma _{\mu ,M-1}\}\colon H_p^{\mu (s)}(\comega)\to \ol
H_p^{s-\operatorname{Re}m}(\Omega )\times \prod _{0\le j<M}B^{s-\operatorname{Re}\mu -j-1/p}_p(\partial\Omega ).\tag4
$$

\endproclaim

Now follow some applications to fractional powers. Let $a>0$ and let $P_a$ equal the power $A^a$ of a strongly elliptic
second-order differential operator $A$ with $C^\infty $-coefficients on $\Omega
_1$ (a special case is $P_a=(-\Delta )^a$). Then $P_a$ is of
order $2a$, of
type $a$, and has factorization index $a$. Theorems 2 and 3 give
e.g.\ the following results in H\"older spaces (where $\dot
C^t(\comega)$ stands for $\{u\in C^t(\Omega _1)\mid \supp
u\subset \comega\}$):

\proclaim{Theorem 4} Let  $u\in \dot H_p^{a-1/p'+0}(\comega )$ for
some $1<p<\infty $ (this
holds if $u\in e^+L_\infty (\Omega )$ when $a<1$, $u\in
\dot C^{a-1+0}(\comega )$ when
$a\ge 1$). The solutions of
$$
r_\Omega P_au=f\tag5
$$
satisfy for $t\ge 0$:
$$
f\in  C^{t+0} (\comega )\implies  u\in e^+d(x)^a C^{t+a-0
}(\comega )\cap C^{t+2a-0}(\Omega ).\tag6
$$
(For $t=0$, $f\in e^+L_\infty
(\Omega )$ suffices.) A solution exists under a finite dimensional
linear condition on  $f$. Moreover,
$$
f\in  C^{\infty } (\comega )\iff u\in e^+d(x)^a C^{\infty 
}(\comega ),\tag7
$$
with Fredholm solvability.
\endproclaim

This theorem is concerned with the homogeneous Dirichlet problem for $P_a$. We can
moreover treat a nonhomogeneous Dirichlet problem (8):

\proclaim{Theorem 5} Let  $u\in H_p^{(a-1)(s)} (\comega )$ with $s>a-1/p'$. 
The solutions of
$$
r_\Omega P_au=f,\quad \gamma _0d(x)^{1-a}u=\varphi ,\tag8
$$
satisfy
$$
\multline
f\in  C^{t+0} (\comega ), \varphi \in C^{t+a+1}(\partial\Omega
)\implies  \\
 u \in  e^+d(x)^{a-1} C^{t+a+1-0
}(\comega )\cap C^{t+2a-0}(\Omega )+\dot C^{t+2a-0}(\comega).
\endmultline\tag9
$$
(For $t=0$, $f\in e^+L_\infty
(\Omega )$ suffices.) A solution exists under a finite dimensional
linear condition on  $\{f,\varphi \}$. Moreover,
$$
f\in  C^{\infty } (\comega ), \; \varphi \in C^\infty (\partial\Omega
)\iff u\in e^+ d(x)^{a-1} C^{\infty 
}(\comega ),\tag10
$$
with Fredholm solvability.
\endproclaim

Ros-Oton and Serra have recently shown in \cite{RS14} for (5) with  $P_a=(-\Delta )^a$,
$0<a<1$, that $f\in L_\infty $ implies $u\in d(x)^aC^\alpha $ for an
$\alpha <\min\{a, 1-a\}$ when $\Omega $ is $C^{1,1}$, by potential
theoretic methods. 
Theorem 4 sharpens this result, allows more general operators,
and extends it to higher
regularity, when $\Omega $ is smooth. We are not aware of any 
published precedents to the other theorems given above. One can also
replace the condition in (8) by a Neumann condition $\gamma
_{a-1,1}u=\psi $ or more general conditions.

\medskip
The theory of $\mu $-transmission $\psi $do's presented here provides a missing link
between, on one hand, Boutet de Monvel's theory of boundary value
problems for integer-order $0$-transmission $\psi $do's, and on the
other hand the very general boundary value theories of other
authors. There is a rich literature; let us for example point to the
works of Schulze
and coauthors, see e.g.\ Rempel-Schulze \cite{RS84},
Harutyunyan-Schulze \cite{HS08} and their references, and the works
of Melrose and coauthors, e.g.\ Melrose \cite{M93}, 
Albin and Melrose 
\cite{AM09} and their references.

\medskip

\noindent{\it Outline.} In Section 1, the relevant function spaces are introduced, including
H\"ormander's $\mu $-spaces,  along with
important order-reducing operators. Section 2 defines the $\mu
$-trans\-mis\-si\-on property and the corresponding boundary behavior for
smooth functions. Section 3 recalls the result of Vishik and Eskin. In
Section 4 we show the Sobolev mapping properties of $\mu
$-trans\-mis\-si\-on operators and deduce the regularity results for
solutions of elliptic homogeneous boundary problems. Section 5 defines
the appropriate boundary operators, and analyzes the structure of the
solution spaces. In Section 6, solvability of
nonhomogeneous elliptic boundary problems is established, with a description 
of parametrices. Finally in Section 7,
consequences are drawn for fractional powers of strongly elliptic differential
operators, and their solvability properties in H\"older spaces.

\head 1. Function spaces \endhead

\subhead 1.1 $L_p$-Sobolev spaces\endsubhead
The function spaces used in \cite{H65} are $L_2$-Sobolev spaces and
their anisotropic variants as introduced in \cite{H63}, together
with a hitherto unpublished interesting case describing a special boundary behavior adapted
to symbols with the $\mu $-transmission property.

In the present paper we generalize this to $L_p$-Sobolev spaces,
mainly of Bessel-potential type, $1<p<\infty $, to which the results
of Eskin's book \cite{E81} were extended in \cite{S94} and \cite{CD01}.
The notation will be a compromise between the nowadays common style
where the regularity exponent $s$ is an upper index without
parentheses, giving room for $p$ as a lower index
(in \cite{H63, H65, H85}, a lower index $(s)$ is used), and on the other
hand H\"ormander's notation of indicating by $\ol H(\rnp)$ resp.\ $\dot
H(\crnp)$ the distributions {\it restricted from} $\R^n$ resp.\ {\it
supported in} $\crnp$. The spaces are all Banach spaces with the
indicated norms.

In the Euclidean space $\R^n$, the points are written
$x=\{x_1,\dots,x_n\}=\{x',x_n\}$, $\rnpm=\{x\mid x_n\gtrless 0\}$,
$\ang x=(1+|x|^2)^\frac12$, and we denote by $[\xi ]$ a smoothed
version of $|\xi |$:$$
[\xi ]\in C^\infty ({\Bbb R}^n,\rp),\; [\xi ]=|\xi |\text{ for }|\xi |\ge 1,\;
[\xi ]\ge \tfrac12\text{ for all }\xi .\tag1.1$$
Restriction from $\R^n$ to $\rnpm$ is denoted $r^\pm$,
 extension by zero from $\rnpm$ to $\R^n$ is denoted $e^\pm$. 

$\F$ denotes the Fourier transformation
$$
(\F f)(\xi )=\hat f(\xi )=\int_{{\Bbb R}^n}e^{-ix\cdot \xi }f(x)\, dx,
$$
 defined on the Schwartz space $\SD({\Bbb R}^n)$ of rapidly decreasing
 $C^\infty $-functions, and extended to distribution in $\SD'(\R^n)$ and
 in function spaces in a well-known way. Note the minus-sign, standard
 in the Western literature, whereas there is usually a
 plus-sign in the definition used in the literature originating from
 Russian and other East-european authors. 

We shall consider classical pseudodifferential operators ($\psi $do's) $P$ of order $m
\in{\Bbb C}$; this means that the symbol has an expansion  in
homogeneous terms $p(x,\xi )\sim \sum_0^\infty p_j(x,\xi )$, where
$p_j$ is homogeneous of degree $m -j$ in $\xi $:
$$
p_{j}(x,t\xi )=t^{m-j}p_j(x,\xi )=t^{\operatorname{Re}m-j}e^{i
\operatorname{Im}m\log t}p_j(x,\xi ),\text{ for }t>0.
$$
 (We just take
one-step polyhomogeneous symbols here, although \cite{H65} allows general
order sequences $m_j$ with $\operatorname{Re}m_j\to -\infty $.) The
operator is defined by
$$
Pu=p(x,D)u=\operatorname{OP}(p(x,\xi ))u
=(2\pi )^{-n}\int e^{ix\cdot\xi
}p(x,\xi )\hat u\, d\xi,\tag1.2 
$$ 
suitably interpreted. Some boundary problems are treated
e.g.\ in \cite{B71, G90, G96, G09}. By truncation to $\rnpm$, $P$ defines
$P_\pm=r^\pm Pe^\pm$.

For $s,t\in \R$ and $1<p<\infty $, the Bessel-potential spaces over
$\R^n$
are defined by
$$
\aligned
H^s_p(\R^n)&=\{u\in \SD'({\Bbb R}^n)\mid \F^{-1}(\ang{\xi }^s\hat u)\in
L_p(\R^n)\},\\
&
\text{ with norm }\|u\|_{H^s_p(\R^n)}=\|u\|_s=\|\F^{-1}(\ang{\xi }^s\hat u)\|_{ L_p(\R^n)},\\
H^{s,t}_p(\R^n)&=\{u\in \SD'({\Bbb R}^n)\mid \F^{-1}(\ang{\xi }^s\ang{\xi '}^t\hat u)\in L_p(\R^n)\},\\
&\text{ with norm }\|u\|_{H^{s,t}_p(\R^n)}=\|u\|_{s,t}=\|\F^{-1}(\ang{\xi }^s\ang{\xi '}^t\hat u)\|_{ L_p(\R^n)}.
\endaligned\tag1.3
$$

The latter anisotropic spaces are used in \cite{H63, G96, G09, CD01};
\cite{S94} includes other anisotropic cases. Note that $H_p^s=H_p^{s,0}$,
and that $H^0_p=L_p$.

The pseudodifferential symbols  $p(x,\xi )$ of order $m\in{\Bbb C}$ are in
$S^{\operatorname{Re}m }_{1,0}(\R^n\times\R^n)$, hence the operators
are continuous from $H^s_p(\R^n)$ to $H^{s-\operatorname{Re}m }_p(\R^n)$ for all
$s\in{\Bbb R}$, as accounted for e.g.\ in \cite{G90}. The continuity
extends to the map from $H^{s,t}_p(\R^n)$ to $H^{s-\operatorname{Re}\mu ,t}_p(\R^n)$ for all $t\in{\Bbb
R}$, cf.\ e.g.\ \cite{CD01}. The operators we consider in this paper
are scalar.

From the spaces in (1.3) we define with a notation extended from
\cite{H63, H65, H85}:
$$
\aligned
\dot H^{s,t}_p(\crnp)&=\{u\in H^{s,t}_p({\Bbb R}^n)\mid \supp u\subset
\crnp \},\\
\ol H^{s,t}_p(\rnp)&=\{u\in \D'(\rnp)\mid u=r^+U \text{ for some }U\in
H^{s,t}_p(\R^n)\},
\endaligned \tag1.4
$$
the first space is a closed subspace of $H^{s,t}_p(\R^n)$, and in
the second space, homeomorphic to $H^{s,t}_p(\R^n)/\dot H^{s,t}_p(\crnm)$, the norm
$$
\|u\|_{\ol H^{s,t}_p(\rnp)} =\inf\{\|U\|_{H^{s,t}_p(\R^n)}\mid
u=r^+U\},\text{ also denoted }\|u\|_{s,t},
$$
is used. $\dot H$ was denoted $\oc H$ in the book \cite{H63} and in \cite{H65}. In
some other texts it is marked as $H_0$ (e.g.\ in \cite{G90}), 
or  $\widetilde H$ (e.g.\ in \cite{E81, T95, S94, CD01}). When $s-1/p$ is integer, Triebel's use of
$\oc H$ in \cite{T95} (first edition 1978) differs from H\"ormander's original
1963 definition. 

The use of both $\ol H$ and $\dot H$ is
practical, since it allows leaving out the indication of the domain $\rnp$. We
recall that $\dot H^{s,t}_p(\crnp)$ and $\ol H^{\,-s,-t}_{p'}(\rnp)$
($1/p'=1-1/p$) are
dual spaces to one another with respect to an extension of the sesquilinear form
$(u,v)=\int_{\rnp}u(x)\overline v(x)\, dx$. 

We shall denote
$$
\bigcup_{\varepsilon >0}\dot H_p^{s+\varepsilon }=\dot
H_p^{s+ 0},\;
\bigcap_{\varepsilon >0}\dot H_p^{s-\varepsilon }=\dot
H_p^{s- 0},\;
\bigcup_{\varepsilon >0}\ol
H_p^{s+\varepsilon }=\ol H_p^{s+ 0},\;
\bigcap_{\varepsilon >0}\ol
H_p^{s-\varepsilon }=\ol H_p^{s-0}.
\tag1.5
$$

The notation $\dot{\SD}(\crnp)$, $\dot{\SD}'(\crnp)$, will be used for
Schwartz functions resp.\ distributions supported in $\crnp$, and 
$\ol{\SD}(\rnp)$, $\ol{\SD}'(\rnp)$, will be used for
Schwartz functions resp.\ distributions restricted to $\rnp$.
Here $\dot{\SD}(\crnp)$ (and $C_0^\infty (\rnp)$) is dense in the spaces $\dot H^{s,t}_p(\crnp)$,
and $\ol{\SD}(\crnp)$ is dense in $\ol H^{s,t}_p(\rnp)$. 

We shall also need the
Besov spaces 
 $B^s_{p}(\R^n)$, which enter as 
range spaces for trace maps, recalling
that for $0<s<2$, 
$$
f\in
B^s_p({\Bbb R}^{n})\iff \|f\|_{L_p}^p+ \int_{\R^{2n}}\frac{|f(x)+f(y)-2f((x+y)/2)|^p}{|x+y|^{n+ps}}\,dxdy<\infty ;
$$
and $B^{s-t}_p(\R^n)=(1-\Delta )^{t/2}B^s_p({\Bbb R}^n)$ for all $t\in{\Bbb R}$.

Embedding, interpolation and other properties are found e.g.\ in Triebel \cite{T95}.

Let $\gamma _j$ denote the trace operator $\gamma _j\colon
u(x',x_n)\mapsto D_n^ju(x',0)$, defined to begin with on smooth
functions: it extends to a continuous linear map
$\gamma _j\colon \ol H^s_p(\rnp)\to B^{s-1/p}_p(\R^{n-1})$, for
$s>1/p$. It is surjective with a continuous right inverse. In fact,
defining the column vector 
$\varrho _M=\{\gamma _0,\dots,\gamma _{M-1}\}$ for a positive integer $M$, we have that 
$$
\varrho _M\colon \ol H^s_p(\rnp)\to\prod_{0\le j<M}B^{s-j-1/p}_p(\R^{n-1})
\text{ for }s>M-1/p,\tag1.6
$$
continous and surjective, having a right inverse (row vector) $\Cal
K_M=\{K_0,\dots, K_{M-1}\}$ (a Poisson
operator, cf.\ \cite{G90}), that in addition is continuous from
$\prod_{0\le j<M}B^{t-j-1/p}_p(\R^{n-1})$ to $\ol H^t_p(\rnp)$ for all
$t\in \R$. As $\Cal K_M$ one can for example take the Poisson operator
$\varphi \mapsto u$ solving
the Dirichlet problem for $(1-\Delta )^M$,$$
(1-\Delta )^Mu=0\text{ in }\rnp,\; \varrho _{M}u=\varphi \text{ on }{\Bbb R}^{n-1}
$$ (an elementary treatment of the case $M=1$ is found in \cite{G09}, Ch.\ 9). 
We shall here use the closely related choice, cf.\ (1.1) ($e^+$ is sometimes left out):
$$
\aligned
&\Cal K_{M}=\{K_0,\dots, K_{M-1}\},\text{ with}\\
K_j\colon \varphi _j &\mapsto \tfrac{(-1)^j}{j!}\F^{-1}_{\xi \to x}\bigl(\hat\varphi _j(\xi ')\partial_{\xi _n}^j ([\xi ']+i\xi _n)^{-1}\bigr)
=\tfrac{i^j}{j!}x_n^j\F^{-1}_{\xi '\to x'}\bigl(e^+r^+e^{-[\xi ']x_n}
\hat\varphi _j(\xi ')\bigr).
\endaligned\tag1.7
$$
It can also be convenient to use (1.7) with $[\xi ']$ replaced by
$\ang{\xi '}$, more closely related to $1-\Delta $.
Still another choice
is given in \cite{H63}, Th.\ 2.5.7 (also recalled in \cite{G96, G09}). 

It is known that there are natural identifications
$$
\aligned
\dot H^s_p(\crnp)&= \{u\in \ol H^s_p(\rnp)\mid \varrho _Mu=0\}, \text{
for }M+1/p>s>M+1/p-1;\\
\dot H^s_p(\crnp)&=  \ol H^s_p(\rnp), \text{
for }1/p>s>1/p-1=-1/p'.
\endaligned \tag1.8
$$

In the borderline case 
$s=1/p$, $\ol H^{1/p}_p(\rnp)$ is strictly larger than $\dot
H^{1/p}_p(\crnp)$; the latter carries the norm $\|u\|_{\ol
H^s_p}+\|x_n^{-1/p}u\|_{L_p}$. However, $C_0^\infty (\rnp)$ is dense
in both of these spaces. (Cf.\ \cite{G90} (2.15)ff.\ and its references.)

The definitions carry over to the manifold situation by
use of local coordinates.

\subhead 1.2 Order-reducing operators\endsubhead
Homeomorphisms between the various spaces  play an important role in
the theory. The operator $\operatorname{OP}(\ang\xi ^\mu )$ defines
homeomorphisms from $H^s_p(\R^n)$ to $H^{s-\operatorname{Re}\mu
}_p(\R^n)$ for all $s\in\R$. Likewise for any $\mu \in{\Bbb C}$, cf.\ (1.1),
$$
\aligned
\Xi ^\mu& =\operatorname{OP}(\chi ^\mu ), \text{ where } \chi ^\mu
=[\xi ]^\mu ,\text{ defines homeomorphisms}\\
\Xi ^\mu &\colon H^s_p(\R^n)\simto H^{s-\operatorname{Re}\mu }_p(\R^n),
\text{ all }s\in\R, \text{ with inverse }\Xi ^{-\mu }. 
\endaligned\tag1.9
$$
In the following, we can either use  $\ang{\xi }, \ang{\xi '}$ as in
\cite{H65}, or replace them by $[\xi ], [\xi ']$ to profit from the
homogeneity. The operators defined by the two choices have the same
mapping properties. The explicit formulas in the following will be
written with $[\xi ']$, since this is useful in the definition of
$\Lambda _\pm^\mu $ further below.

For the spaces defined relative to
$\rnpm$, there are several interesting choices. One is the simple family$$
\chi ^\mu _+=([\xi ']+i\xi _n)^\mu ,\text{ resp.\ }\chi ^\mu _-
=([\xi ']-i\xi _n)^\mu ,\quad \operatorname{OP}(([\xi ']\pm i\xi _n)^\mu )=\Xi
^\mu _{\pm}, \tag1.10
$$   
(or, if needed, the corresponding formulas with $\ang{\xi '}$). Here
$\chi ^\mu _+$ (resp.\ $\chi ^\mu _-$) extends analytically as a
function of $\xi _n$ into
${\Bbb C}_-=\{\operatorname{Im}\xi _n<0\}$ resp.\ ${\Bbb C}_+=\{\operatorname{Im}\xi
_n>0\}$. (The imaginary halfspaces play the opposite roles in the works
\cite{E81, S94, CD01} because of the opposite sign in the definition
of $\F$.) Since $\chi _+^\mu $ extends
analytically to
$\operatorname{Im}\xi _n<0$, the operator $\Xi ^\mu _+$ preserves support in
$\crnp$; hence
we have for all $s\in\R$ that $$
\Xi^\mu _+\colon \dot H^s_p(\crnp)\simto
\dot H^{s-\operatorname{Re}\mu }_p(\crnp),\text{ with inverse }\Xi
^{-\mu }_+.\tag1.11
$$
The adjoint mapping  is $\Xi^{\overline\mu }_{-,+}
\colon \ol H^{\,-s+\operatorname{Re\mu }}_{p'}(\rnp)\simto
\ol H^{\,-s }_{p'}(\rnp)$; this shows for general $s,p,\mu $:
$$
\Xi ^\mu _{-,+}
\colon \ol H^s_p(\rnp)\simto
\ol H^{s-\operatorname{Re}\mu }_p(\rnp),\text{ with inverse }\Xi ^{-\mu }_{-,+}.
\tag1.12
$$

\example{Remark 1.1} For $s>-1/p'$, $\Xi ^\mu _{-,+}$ in (1.12) identifies with
$r^+\Xi^\mu _-e^+$ ($e^+$ is only defined then). For lower
$s$, the mapping in (1.12) can be understood, besides being a
specific adjoint,  as  the extension by
continuity from the operator defined on the dense subspace $\ol{\SD}(\rnp)$ (as noted in
\cite{GK93}, p.\ 174).
There is also a third formulation worth mentioning, used in
\cite{E81}, namely that for any extension operator
$\ell:\ol H^s_p(\rnp)\to H^s_p({\Bbb R}^n)$ with $r^+\ell=\operatorname{Id}$, 
$$
\Xi ^\mu _{-,+}f=r^+\Xi ^\mu _{-}\ell f.\tag1.13
$$
This holds since $r^+\Xi ^\mu _{-}g=0$ for any distribution $g$ supported
in $\crnm$, using that since $\chi ^\mu _-$ extends analytically to
$\operatorname{Im}\xi _n>0$, the operator $\Xi ^\mu _{-}$ preserves support in
$\crnm$. The formula (1.13) is independent of the choice of $\ell$.
\endexample

The symbols $\chi ^\mu _\pm$ are not truly pseudodifferential
(although the $\operatorname{OP}(\chi ^\mu _\pm)$ have a good meaning by Lizorkin's
criterion, cf.\ e.g.\ \cite{G90}), since the higher
$\xi '$-derivatives do not have the correct fall-off for $|\xi
|\to\infty $. But there exists another choice with true $\psi $do symbols
given  in \cite{G90} (inspired from the unpublished \cite{F86}), that also has the above
mapping properties. 
Define
$$
\lambda _{\pm}^\mu =(\lambda _{\pm}^1)^\mu ,\quad \lambda _{-}^1 =
[\xi ']\psi \Bigl(\frac{\xi _n}{a[\xi ']}\Bigr)- i\xi
_n,\quad \lambda _+^1=\ol {\lambda _-^1},
\tag 1.14
$$
with $\psi \in\SD(\R)$ having $\psi (0)=1$ and $\supp\F^{-1}\psi \subset
\crm$. We set $\psi (\pm\infty )=0$, then $\psi $ is $C^\infty $ on
the extended real axis. Here the constant $a>0$ is chosen so large that the negative powers
are well-defined, cf.\ \cite{G90} pp.\ 317-322. 
The functions
$\lambda ^\mu _+$ (resp.\ $\lambda ^\mu _-$) extends analytically into
$\{\operatorname{Im}\xi _n<0\}$ resp.\ $\{\operatorname{Im}\xi
_n>0\}$. 
 Denoting
$\operatorname{OP}(\lambda ^\mu _\pm)=\Lambda 
^\mu _{\pm}$, we have for all $s\in{\Bbb R}$ that $$
\aligned
\Lambda^\mu _+&\colon \dot H^s_p(\crnp)\simto
\dot H^{s-\operatorname{Re}\mu }_p(\crnp),\text{ with inverse }\Lambda
^{-\mu }_+,\\
\Lambda ^\mu _{-,+}&\colon \ol H^s_p(\rnp)\simto
\ol H^{s-\operatorname{Re}\mu }_p(\rnp),\text{ with inverse }\Lambda ^{-\mu }_{-,+};
\endaligned
\tag1.15
$$
here $\Lambda ^\mu _{-,+}$ is the adjoint of $\Lambda^{\overline\mu }_+\colon \dot H^{-s+\operatorname{Re}\mu }_{p'}(\crnp)\simto
\dot H^{-s }_{p'}(\crnp)$, and again there are interpretations as in Remark 1.1. The proofs are given in \cite{G90},  (cf.\ (4.11), (4.24)
there) using that 
 for $a$ taken sufficiently
large in (1.14) (as we assume),
$$
\eta ^\mu _{\pm} (\xi )=(\lambda ^1 _{\pm}(\xi )/\chi ^1_{\pm}(\xi ))^\mu  
=1+q^\mu _{\pm}(\xi )\text{ with } |q^\mu _{\pm}(\xi )|\le \tfrac12,\tag1.16
$$
analytic for  $\operatorname{Im}\xi _n\lessgtr 0$; they define  $\psi
$do's  
$\eta ^\mu _{\pm} (\xi ',D_n)=\operatorname{OP}_n(\eta ^\mu _{\pm} (\xi ',\xi _n))$ of
order 0  that are homeomorphisms in $L_2({\Bbb R})$, uniformly in $\xi
'$. Since  they preserve support in $\crpm$ respectively
(and the inverses do so too), $r^{\pm}\eta ^\mu _{\pm} (\xi ',D_n)e^{\pm}$ are
homeomorphism in
$L_2(\rpm)$, respectively. This allows transferring the mapping properties of the
$\Xi^\mu _\pm$ to the $\Lambda ^\mu _{\pm}$, cf.\ \cite{G90}. The
operators 
$\Xi ^\mu
_{+}$, $\Lambda ^\mu _{+}$ and $\eta ^\mu _{+} (\xi ',D_n)$ belong to
the so-called ``plus-operators'' of
Eskin \cite{E81}, and the operators
$\Xi ^\mu
_{-}$, $\Lambda ^\mu _{-}$ and $\eta ^\mu _{-} (\xi ',D_n)$ belong to
the  ``minus-operators''. The symbols are said to be ``plus-symbols''
resp.\ ``minus-symbols''. (The sub-indices $\pm$ here should not
be confounded with the $\pm$ used to indicate truncation --- added on
as an extra index.)

In addition to what was shown in \cite{G90}, we observe:

\proclaim{Lemma 1.2} Let $Y^\mu _{+}=\OP(\eta ^\mu _+(\xi ))$, then $Y^\mu _{+,+}=r^+Y^\mu _+e^+$ is a
homeomorphism of $\ol H^{s,t}_p(\rnp)$ onto itself for all
$s,t\in{\Bbb R}$.
For any $s,t\in{\Bbb R}$,
$$
\|r^+\Xi ^\mu _+u\|^{\;}_{H^{s,t }_p(\rnp)}\simeq \|r^+\Lambda ^\mu
_+u\|^{\;}_{H^{s,t }_p(\rnp)}
.\tag1.17
$$
 The equivalence also holds if $[\xi ']$ is replaced by $\ang{\xi '}$
 in the definition of $\Xi ^\mu _{+}$.
\endproclaim

\demo{Proof} The proof needs some care, because $Y^\mu _+$ is not a
standard $\psi $do on $\R^n$; however it is so at the one-dimensional level
where we just use the definition with respect to $\xi _n$. Here the
Boutet de Monvel calculus on ${\Bbb R}$ shows that $r^+\eta ^\mu _+(\xi ',D
_n)e^+$ is a homeomorphism in $\ol H^m_2(\rp)$ with inverse
$r^+\OP_n((\eta ^\mu _+(\xi ))^{-1})e^+$ for all $m\in
{\Bbb Z}$, since the left-over operators such as $G^+(\OP_n(\eta ^\mu _+))G^-(\OP_n(\eta  _+^{-\mu }))$
arising in the composition have the $G^-$-factor equal to 0, hence vanish. The norms are
bounded in $\xi '$. Interpolation extends the homeomorphism property
 to all real $s$.

Estimating the norms simply by Fourier transformation, we find for
$p=2$ that the full operator $r^+Y^\mu _+e^+$ is a homeomorphism in $\ol
H^{s,t}_2(\rnp)$ with inverse $r^+(Y^\mu _+)^{-1}e^+$. Both $Y^\mu _+$ and $(Y^\mu _+)^{-1}=Y^{-\mu }_+$ are
continuous in $H^{s,t}_p(\R^n)$ by Lizorkin's criterion. The
$L_2$-calculations apply in particlar to functions $u\in
\overline{\SD}(\rnp)$, showing that $r^+Y^\mu _+e^+u=r^+e^+Y^\mu _+u$,
$r^+Y^{-\mu } _+e^+u=r^+e^+Y^{-\mu }_+u$ for such $u$; this extends to $u\in
\ol H^{s,t}_p(\rnp)$ by closure, and completes the proof of the
homeomorphism property.

Now
$$
r^+\Lambda _+^\mu u= r^+Y^\mu _+\Xi _+^\mu  u=r^+Y^\mu _+e^+r^+\Xi _+^\mu u,
$$
where the corresponding term with $e^-r^-$ in the middle vanishes
since 
$r^-\Xi _+^\mu u$ does
so. 
Then in view of the homeomorphism property of $r^+Y^\mu _+e^+$,
$$
\|\Lambda
_+^\mu u\|_{s,t} \le C
\|\Xi
_+^\mu u\|_{s,t} ,
$$
Similarly, an inequality the other way follows by use of $Y^{-\mu }_+$.

For the last statement, the operators $\operatorname{OP}(([\xi ']+i\xi
_n)^\mu )$ and $\operatorname{OP}((\ang{\xi '}+i\xi _n)^\mu )$ can be
compared in a similar way, since $\bigl(([\xi ']+i\xi _n)/(\ang{\xi
'}+i\xi _n)\bigr)^\mu = \bigl(1+([\xi ']-\ang{\xi '})/(\ang{\xi
'}+i\xi _n)\bigr)^\mu $ is an invertible plus-symbol of order 0.
\qed
\enddemo

It is important to observe that the operators $\Lambda ^m _+$,
$m\in{\Bbb Z}$, that act homeomorphically in the scale $\dot H_p^s(\crnp)$,
can also be applied to the scale $\ol H_p^s(\rnp)$ for $s>-1/p'$ after truncation, $\Lambda
^m_{+,+}=r^+\Lambda ^m_+e^+$, since they belong
to the Boutet de Monvel calculus.  But here they must in general be
supplied with trace or Poisson operators to define
homeomorphisms. E.g.\ for integer $m>0$,
$$
 \pmatrix \Lambda ^m_{+,+}\\ \quad\\ \varrho _m
\endpmatrix \colon \ol H^s_p(\rnp) \simto \matrix 
\ol H^{s-m}_p(\rnp)\\ \times  \\ \prod _{0\le j<m} B^{s-j-1/p}_p(\R^{n-1})
\endmatrix, \text{ when }s>m-1/p'\tag1.18
$$
(shown in \cite{G90}, Th.\ 4.3); it is an elliptic boundary value problem. A similar mapping property holds with $\Xi ^m_+$ instead
of $\Lambda ^m_+$.

\medskip

The construction of these operators extends to the manifold situation, by the
method described in \cite{G90}.
Let $\comega $ be a compact $n$-dimensional $C^\infty $ manifold with
interior $\Omega $ and boundary $\partial\Omega =\Sigma $, and let $E$
be a Hermitean $C^\infty $ vector bundle over $\comega$ of dimension $N$, its restriction to
$\Sigma $ denoted $E'$. We can assume that $\comega$ is
smoothly embedded in a compact boundaryless $n$-dimensional manifold
$\Omega _1$ (e.g.\ the double of $\comega$) such that $\Sigma $
is the boundary of $\Omega $ there, and we assume that $E$ is the restriction to
$\comega$ of a smooth vectorbundle $E_1$ given over $\Omega _1$. Then there
is a standard way to generalize the definitions of Sobolev spaces over
$\R^n$, $\rnpm$, to
spaces of distributions over $\comega $, $\Sigma $, $\Omega _1$,
valued in the bundles, by use of local trivializations. The definition
of $\psi $do's likewise generalizes to the manifold and vector bundle
situation. In the present paper, our application deals with
scalar  $\psi $do's, so we shall drop the vector bundle aspect to
simplify notations, but declare at this point that the constructions
of order-reducing operators generalize to bundles as in \cite{G90},
easily taken up when needed. We denote by $r_\Omega $, or for brevity
$r^+$, the restriction
from $\Omega _1$ to $\Omega $, and by $e_\Omega $ or $e^+$ the extension from
$\Omega $ by zero on $\Omega _1\setminus\comega$. For an operator $P$
over $\Omega _1$, we denote $r_\Omega Pe_\Omega $ (also called
$e^+Pr^+$) by $P_\Omega $ or $P_+$.

\proclaim{Theorem 1.3} There exists a family of elliptic $\psi $do's $\Lambda
_+^{(\mu )}$on $\Omega _1$, classical of order $\mu $ and with
principal symbol $\lambda ^\mu _+$ at the boundary of $\Omega $,  preserving support in $\comega$ and defining homeomorphisms
$$
\Lambda^{(\mu )}_+\colon \dot H^s_p(\comega )\simto
\dot H^{s-\operatorname{Re}\mu }_p(\comega ),\tag1.19
$$
for all $s\in{\Bbb R}$, with inverses $(\Lambda^{(\mu )}_+)^{-1}$
likewise preserving support in $\comega$. The family of adjoints 
are classical elliptic operators $\Lambda ^{(\overline \mu )}_-$, with
principal symbol $\lambda _-^{\overline\mu }$ at the boundary of $\comega$,
such that $\Lambda ^{(\mu )}_{-,+}=r^+\Lambda ^{(\mu )}_-e^+$ are homeomorphisms
$$
\Lambda ^{(\mu )}_{-,+}\colon \ol H^s_p(\Omega )\simto
\ol H^{s-\operatorname{Re}\mu }_p(\Omega ),\tag1.20
$$
for all $s\in {\Bbb R}$, with inverses $((\Lambda ^{(\mu )}_-)^{-1})_+$. 
\endproclaim

\demo{Proof} The construction is explained in detail in \cite{G90},
Sections 4 and 5, which we use with minor adaptations that we shall
explain here. We provide $\Omega _1$ and $\Sigma $ with
Riemannian metrics, such that a tubular neighborhood $\Omega _2$ of
$\Sigma $ in $\Omega _1$ is isometric with $\Sigma \times \,]-2,2[\,$;
the coordinates in $\Sigma $ resp.\ $\,]-2,2[\,$ will be denoted $x'$
and $x_n$, and we write $\Sigma _c=\Sigma \times\,]-c,c[\,$ for $c\le
2$.  Fix $\mu $. In the definition of $\lambda ^\mu _\pm$ (1.14) we can insert an
extra parameter $\zeta \ge 0$ (called $\mu $ in \cite{G90}), defining
$$
\lambda _{\pm,\zeta }^\mu =(\lambda _{\pm,\zeta }^1)^\mu ,\quad
\quad \lambda _{-,\zeta }^1 =
[(\xi ',\zeta )]\psi (\xi _n/a([(\xi ',\zeta )])- i\xi
_n,\quad \lambda _{+,\zeta }^1=\ol {\lambda _{-,\zeta }^1}.
\tag 1.21
$$
Now the construction of the $\psi $do $\Lambda ^{(\mu
)}_{+,\zeta }$ defined on $\Omega _1$ is carried out similarly to
the description in
\cite{G90} around (5.1), 
 using
$\lambda _{+,\zeta }^\mu $ near the boundary and
$[(\xi ,
\zeta )]^{\mu }$ at a distance from the boundary:
$$
\lambda ^{(\mu )}_{+,\zeta }=(\lambda _{+,\zeta }^1)^{\mu \alpha (x_n)}[(\xi ,\zeta )]^{\mu (1-\alpha (x_n))} 
$$
on $\Sigma _2$, extended by $[(\xi ,\zeta )]^{\mu } $ on the rest of
$\Omega _1$; here $\alpha (x_n)\in C^\infty ({\Bbb R},[0,1])$ equal to
$1$ on $[-1,1]$ and 0 on the complement of $[-\frac32,\frac32]$. The
symbol extends analytically to $\operatorname{Im}\xi _n<0$. The
operator $\Lambda ^{(\mu )}_{+,\zeta }$ is pieced together from this by use of a
finite partition of unity subordinate to a covering of $\Omega _1$ by
open sets in $\Sigma _{\frac34}$ and open sets in $\Omega _1\setminus
\Sigma _{\frac12}$, whereby $\Lambda ^{(\mu )}_{+,\zeta }$ preserves
support in $\comega$.

The construction with $\mu $ replaced by $-\mu $ gives the operator
$\Lambda ^{(-\mu )}_{+,\zeta }$, likewise elliptic on $\Omega _1$ and
preserving support in $\comega$. Now
$$
\Lambda ^{(\mu )}_{+,\zeta }\Lambda ^{(-\mu )}_{+,\zeta }=I+U_1(\zeta ),\quad\Lambda ^{(-\mu )}_{+,\zeta }\Lambda ^{(\mu )}_{+,\zeta }=I+U_2(\zeta ),\tag1.22
$$
with $U_1,U_2$ of order $-1$, hence compact operators in $H^t_p(\Omega
_1)$ for all $t,p$; they also preserve support in
$\comega$. Standard elliptic theory shows that $\Lambda ^{(\mu
)}_{+,\zeta }$ is a Fredholm operator from $H^s_p(\Omega _1)$ to
$H^{s-\operatorname{Re}\mu }_p(\Omega _1)$ for all $s,p$, with a
finite dimensional $C^\infty $ kernel and range complement independent
of $s,p$. We have in particular that $\Lambda ^{(\mu )}_{+,\zeta
}$ maps $\dot H^s_p(\comega)$ into $\dot H^{s-\operatorname{Re}\mu }_p(\comega)$, 
and $\Lambda ^{(-\mu )}_{+,\zeta }$ maps
the other way, with (1.22) valid there, so $\Lambda ^{(\mu
)}_{+,\zeta}$ is Fredholm between those spaces, with a finite
dimensional $C^\infty $ kernel $K_1$ and range complement $K_2$
independent of $s,p$.  The idea with the parameter $\zeta $ is that we can
apply the calculus of \cite{G96} (just for $\psi $do symbols), where our symbols are of {\it
regularity} $\nu =+\infty $ as functions of $(\xi ,\zeta  )$; then the
norms of $U_1$ and $U_2$ are $\le \frac12$ for $\zeta $ sufficiently
large, so that $I+U_1$ and $I+U_2$ are invertible, and it follows that
$\Lambda ^{(\mu )}_{+,\zeta }$ over $\comega$ is invertible for large $\zeta $. Since
it depends continuously on $\zeta $, it follows that $\Lambda ^{(\mu
)}_{+,0}$ has index 0. For $p=2$, the kernel and range complement are
spanned by orthonormal systems of smooth functions $\{\varphi
_1,\dots,\varphi _N\}$   and $\{\psi _1,\dots,\psi _N\}$ supported in $\comega$, and when we
define the order $-\infty $ operator $\Psi $ by  $\Psi u=\sum_{j,k=1}^N\psi
_j(u,\varphi _k)$, 
$$
\Lambda ^{(\mu )}_+=\Lambda ^{(\mu )}_{+,0}+\Psi ,
$$
  has the desired bijectiveness property. 

An operator $\Lambda ^{(\mu )}_{-,+}$ with the desired properties is
now found as the adjoint of $\Lambda ^{(\overline\mu )}_+$ in (1.19),
in the same way as for $\rnp$.\qed
\enddemo

For negative $s$ in (1.20) the operator is understood as in Remark 1.1.
The assertion (1.18) generalizes to these operators. More properties
are shown in Example {\rm 2.8} later. 

It is the introduction of these
$\psi $do's that allows a relatively elegant deduction of solvability
properties for the equations we consider in this paper. They had not
been found when \cite{H65} was written (and there is a remark there
that such operators would be helpful).
\medskip

Occasionally we shall refer to the spaces $C^t(\comega )$ and
$C^t(\Omega )$ for $t\ge 0$; in integer cases they are the usual
spaces of functions with continuous derivatives up to 
order $t$ on $\comega$ resp.\ $\Omega $, and when $t=k+s$, $k\in{\Bbb
N}_0$, $s\in \,]0,1[\,$, they are the H\"older spaces also denoted
$C^{k,s}(\comega)$ resp.\ $C^{k,s}(\Omega)$. We denote 
$\bigcup_{\varepsilon >0}C^{t+\varepsilon }=C^{t+0}$, and
$\bigcap_{\varepsilon >0}C^{t-\varepsilon }=C^{t-0}$ if $t>0$. There are embeddings
$$
\ol H^t_p(\Omega )\subset C^{t-n/p-0}(\comega)\text{ when }t>n/p,\quad C^{t+0}(\comega)\subset
\ol H^{t}_p(\Omega )\text{ when }t\ge 0;\tag1.23
$$
in the first embedding,``$-0$'' can
be left out if $t-n/p$ is not integer, in the second we assume $\comega
$ compact. We shall denote $
\{u\in C^t(\Omega _1)\mid \supp
u\subset \comega\}=\dot C^t(\comega)$.

\subhead 1.3 H\"ormander's $\mu $-spaces \endsubhead
In the notes \cite{H65} there are introduced 
(for $p=2$) the following spaces
that mix the features of the supported and the restricted Sobolev spaces in a 
particular way by use of the mappings $\Xi ^\mu _+$. (Actually,
\cite{H65} uses $(\ang{D'}+\partial_n)^\mu $ instead of
$\Xi _+^\mu =([D']+\partial_n)^\mu $; they are equivalent.) 

\proclaim{ Definition 1.4} Let $\mu \in{\Bbb C}$, and let
$s>\operatorname{Re}\mu -1/p'$. An element $u\in \dot{\SD} '(\crnp)$ is in
$H^{\mu (s)}_p(\crnp)$ if and only if  $\Xi ^\mu _+u\in \dot H_p^{ -1/p'+0}(\crnp)$ and 
$$
\|r^+\Xi ^\mu _+u\|^{\;}_{\ol H_p^{s-\operatorname{Re}\mu }(\rnp)}<\infty ;\tag1.24
$$
the topology is defined by the norm {\rm (1.24)}, also denoted $\|u\|_{\mu (s)}$.

In this definition, $\Xi_+^\mu $ can be replaced by $\Lambda _+^\mu$.
\endproclaim

The last statement is justified by the properties shown in Section 1.2, in
particular Lemma {\rm 1.2}.

The condition $\Xi ^\mu _+u\in \dot H_p^{ -1/p'+0}(\crnp)$ can also be expressed as
$$
u\in \dot H_p^{\operatorname{Re}\mu  -1/p'+0}(\crnp),
$$
in view of the homeomorphism properties (1.11).
 Note that the inequality in (1.24) implies, since $s-\operatorname{Re}\mu >-1/p'$,
that the elements satisfy for  $0<\varepsilon <\min\{1, s-\operatorname{Re}\mu +1/p'\}$:
$$
\Xi ^\mu _+u\in \ol H^{\varepsilon -1/p'}_p(\rnp)\simeq \dot
H^{\varepsilon -1/p'}_p(\crnp),
\tag1.25
$$
using the identification of $r^+v$ and $e^+r^+v$ in spaces with
$-1/p'<s<1/p$, cf.\ (1.8). So the norm (1.24) is stronger than the
norm on the spaces in (1.25), which need not be mentioned in the
definition of the topology.

If $s<\operatorname{Re}\mu +1/p$, the condition in (1.24) reduces to
$\Xi _+^\mu u\in \dot H^{s-\operatorname{Re}\mu }_p(\crnp)$; therefore
$$
H^{\mu
(s)}_p(\crnp)=\dot H^s_p(\crnp)\text{ when }-1/p'<s-\operatorname{Re}\mu <1/p,\tag1.26
$$ 
  and $C_0^\infty (\rnp)$ is dense in the space. 
When $s$ is larger, (1.24) gives a
nontrivial restriction on $u$.

We can then extend the definition to all $s$, consistently with the above:

\proclaim{Definition 1.5} Let $\mu \in\C$, and let
$s<\operatorname{Re}\mu +1/p$. Then we define
$$
H^{\mu
(s)}_p(\crnp)=\dot H^s_p(\crnp).\tag1.27
$$ 
\endproclaim

Note that $\dot H^s_p(\crnp)\subset H_p^{\mu (s)}(\crnp)$ holds for all
$s$ and $\mu $.

\example{Example 1.6} Let $\mu =m\in{\Bbb N}$ and
$s> m-1/p'$. Then $u\in H_p^{m(s)}$ if and only if $u\in \dot H_p^{m-1/p'+0}$
and $r^+([D']+iD_n)^mu\in \ol H_p^{s-m}$. The first condition
implies that  $\varrho _mu=0$, and the second condition holds
if $u\in \ol H_p^{s}$. 
The second condition can also be written $\Lambda ^m_{+,+}u\in \ol
H_p^{s-m}$, and in view of the ellipticity of the system $\{\Lambda
^m_{+,+} \varrho _m\}$ in the Boutet de Monvel calculus, cf.\ 
(1.18), we see that $u$ must lie in $\ol H_p^{s}$.

This shows that $H_p^{m(s)}=\{u\in \ol
H^s_p\mid \varrho _mu=0\}$. Note that for $s>m+1/p$, the space is a
proper subspace of $\ol H_p^s$, different from $\dot H_p^s$.
\endexample

This example is still within the Boutet de Monvel calculus; the
novelty of the spaces $H_p^{\mu (s)}$ lies more in what happens for
noninteger $\mu $.

The following observation will be very useful:

\proclaim{Proposition 1.7} Let $s>\operatorname{Re}\mu -1/p'$. 
The mapping
$r^+\Xi _+^\mu $ is a homeomorphism of $H_p^{\mu (s)}(\crnp)$ onto
$\ol H_p^{s-\operatorname{Re}\mu }(\rnp)$ with inverse $\Xi _+^{-\mu
}e^+$. In particular, $H_p^{\mu (s)}(\crnp)$ is a Banach space.

The analogous result holds with $\Lambda ^\mu _+$-operators, and with
$\Xi ^\mu _+$-operators where $[\xi ']$ is replaced by $\ang{\xi '}$.
\endproclaim

\demo{Proof} By definition, $r^+\Xi _+^\mu $ is continuous.

Surjectiveness is seen as follows: Let $v\in \ol
H_p^{s-\operatorname{Re}\mu }$, and set $w=\Xi _+^{-\mu }e^+v$. Then
$\Xi _+^\mu w=\Xi _+^\mu \Xi _+^{-\mu }e^+v=e^+v$. Since
$s-\operatorname{Re}\mu >-1/p'$, $e^+v\in \dot H_p^{-1/p'+0}$, so $\Xi
^\mu _+w\in  \dot H_p^{-1/p'+0}$ as required in Definition 1.4. Moreover, 
$$
r^+\Xi _+^\mu w=r^+\Xi _+^\mu \Xi _+^{-\mu }e^+v=r^+e^+v=v
$$
is in $\ol H_p^{s-\operatorname{Re}\mu }$ by hypothesis, so $v$ is the
image of $w\in
H_p^{\mu (s)}$.

The injectiveness.
When $u$ satisfies the hypotheses of Definition
{\rm 1.2}, then $u$ is reconstructed from $v=r^+\Xi ^\mu _+u$ as
follows: 
Since $\Xi ^\mu _+u\in \dot H_p^{ -1/p'+0}(\crnp)$, we can write
$$\Xi _+^\mu u = e^+r^+\Xi _+^\mu u+e^-r^-\Xi _+^\mu u.\tag1.28
$$
Here $r^-\Xi _+^\mu u=0$, since  $\Xi ^\mu _+$ preserves support in
$\crnp$. Hence
$$
u=\Xi _+^{-\mu }\Xi _+^\mu u=\Xi _+^{-\mu }e^+r^+\Xi _+^\mu u=\Xi _+^{-\mu }e^+v.
$$

Thus $r^+\Xi _+^\mu $ is an isometry of $H_p^{\mu (s)}$ onto $\ol
H_p^{s-\operatorname{Re}\mu }$, with inverse $\Xi _+^{-\mu }e^+$. In
particular, $H_p^{\mu (s)}$ is a Banach space.

The proof for $\Lambda _+^\mu $ and for the other version of $\Xi
_+^\mu $ goes in the same way.
\qed
\enddemo

The spaces can also be defined in the manifold situation. By use of
the operators $\Lambda ^{(\mu )}_{\pm}$ introduced in Theorem 1.3, we can
formulate the definition as follows:

\proclaim{ Definition 1.8} Let $\mu \in{\Bbb C}$. When
$s>\operatorname{Re}\mu -1/p'$, then 
$H^{\mu (s)}_p(\comega)$ consists of the elements $u\in \dot{\E}'(\comega)$
such that $\Lambda  ^{(\mu )} _+u\in \dot H_p^{ -1/p'+0}(\comega)$ and   
$$
\|r_\Omega \Lambda  ^{(\mu )} _+u\|^{\;}_{\ol H_p^{s-\operatorname{Re}\mu }(\Omega )}<\infty ;\tag1.29
$$
it is a Banach space with the norm {\rm (1.29)}, also denoted $\|u\|_{\mu (s)}$.

When $s<\operatorname{Re}\mu +1/p$, we define
$$
H^{\mu
(s)}_p(\comega)=\dot H^s_p(\comega).\tag1.30
$$ 
\endproclaim

Here the space $\dot{\E}'(\comega)$ denotes the distributions supported
in $\comega$ (compactly supported in $\comega $, if $\comega$ is
allowed to be merely paracompact).
Again we observe that the norm in (1.29) is stronger than the norm in $\dot H^{\varepsilon
-1/p'}_p(\comega)$ for small $\varepsilon $, and that the space equals $\dot
H^s_p(\comega)$ when $-1/p'<s-\operatorname{Re}\mu <1/p$, so that the
last part of the definition allowing lower values of $s$ is consistent
with the first part.
Also Proposition 1.7 extends.

There are of course embeddings $$
H^{\mu (s)}_p\subset H^{\mu (s')}_p\text{
for }s'<s.\tag1.31
$$ On the other hand, embeddings between spaces with
different $\mu ,\mu '$ do not hold in general. An exception is when
$\mu -\mu '$ is integer, see Proposition 4.3 later.

The structure of the spaces will be further described below, 
particularly their importance for $\psi $do's with the transmission
property of type $\mu $.

\example{Remark 1.9}
In \cite{H65}, $H^{\mu (s)}_2(\comega)$ is {\it defined} as the
completion of $\Cal E_\mu (\comega)$ in the topology defined by the
seminorms $u\mapsto \|r^+Pu\|_{\ol H_2^{s-\operatorname{Re}m}(\Omega
)}$, where $P$ runs through the operators of type $\mu $ and any order
$m\in{\Bbb C}$. The proof that this is equivalent with Definition 1.4
(when localized) fills a large section. It is here covered by
Proposition 4.1 and
Theorem 4.2 below.
\endexample

\head 2. The $\mu $-transmission condition  \endhead

The $\mu $-transmission condition is defined and characterized in
\cite{H85} at the end of Section 18.2. Since the explanation is quite
compressed there, we have incorporated some of the original detailed deductions
from \cite{H65} here, slightly modified if necessary. (We remark that the conventions
in \cite{H85} are a little different from here: The
space called $\Cal C_\mu ^\infty $ there  on pp.\ 110--111 is the same
as $\Cal E_{-\mu -1}$ here, and $\mu $ in Th.\ 18.2.18 there
corresponds to $-\mu $ in Definition 2.5 below.)

Let $\Omega _1$ be a fixed paracompact $C^\infty $ manifold, and let $\Omega $
be an open subset of $\Omega _1$ with a $C^\infty $ boundary $\partial\Omega
$. Our purpose is to study boundary problems for the
pseudodifferential operator $P$ in $\Omega $. This means that we
shall look for distributions $u$ with support in $\comega$ such that
$Pu=f$ is given in $\Omega $ and $u$ satisfies some conditions on
$\partial\Omega $ in addition. In particular we shall make a detailed
study of the regularity of $u$ at the boundary when $f$ and the
boundary data are smooth. Examples involving $\alpha $-potentials due
to M.\ Riesz and extended in part by Wallin show that one should not
expect $u$ to be smooth up to the boundary but that one has to expect
$u$ to behave as the distance to the boundary raised to some
power. This leads us to define a family of spaces of distributions
$\Cal E_\mu $ as follows.

\proclaim{Definition 2.1}
If $\operatorname{Re}\mu >-1$ and if $d$ is a real valued function in
$C^\infty (\Omega _1)$ such that 
$$
\Omega =\{x\mid d(x)>0\}\tag2.1
$$
and $d$ vanishes only to the first order on $\partial\Omega $, then
$\Cal E_\mu (\comega )$ consists of all functions $u$ such that $u=0$
in $\complement\comega$ and $u=d^\mu v$ in $\comega$ for some $v\in
C^\infty (\comega)$. 

For lower values of $\operatorname{Re}\mu $, $\Cal E_\mu $ is defined
 successively so that $\Cal
 E_{\mu -1}$ is always the linear hull of the spaces $D\Cal E_\mu $
 when $D$ varies over the first order differential operators with
 $C^\infty $ coefficients. 

\endproclaim

This definition is independent of the choice of
$d$, for if $d_1,d_2$ are two functions with the required properties,
the quotient $d_1/d_2$ is positive and infinitely differentiable. 

To justify the second part of the definition we note that if $D$
is a first order differential operator with $C^\infty $ coefficients,
and if $\operatorname{Re}\mu >0$, then $D\Cal E_\mu \subset \Cal E_{\mu
-1}$, for $D(d^\mu  v)=d^{\mu -1}V$ for some $V\in C^\infty $. The
linear hull of the spaces $D\Cal E_\mu $ when $D$ varies is in fact
equal to
$\Cal E_{\mu -1}$. It is sufficient to prove that it contains any
element in $\Cal E_{\mu -1}$ with support in a coordinate patch where
$\Omega $ is defined by $x_n>0$. Then we can take
$D=\partial/\partial x_n$, noting that if $v\in C^\infty $ then
$$
\int_0^{x_n} t^{\mu -1}v(x', t)\, dt = x_n^\mu V(x),
$$
 where
$$
V(x)=\int_0^1 t^{\mu -1}v(x',x_n t)\, dt 
$$
 is a $C^\infty $ function. If $u=x_n^{\mu -1}v $ and $U=x_n^\mu V\chi
 $, both functions being defined as 0 when $x_n<0$, and $\chi \in
 C_0^\infty $ is 1 in a neighborhood of $\operatorname{supp}u$, then
 $u=\partial U/\partial x_n$ is a $C^\infty $ function on $\rnp$ with
 support in $x_n\ge 0$, so $u\in \partial \Cal E_\mu /\partial x_n+
 \Cal E_\mu $. 
It is thus legitimate to define $\Cal E_\mu $
 successively for decreasing $\operatorname{Re}\mu $ as indicated.

The spaces $\Cal E_\mu $ so obtained have
 the local property that $u\in \Cal E_\mu (\comega)$ and $\varphi \in
 C^\infty (\Omega_1)$ implies that $\varphi u\in \Cal E_\mu (\comega)$. In
 fact, if $D$ again denotes a first order differential operator we have
$$
\varphi D\Cal E_{\mu +1}\subset D\varphi \Cal E_{\mu +1}+\Cal E_{\mu
+1}\subset D\Cal E_{\mu +1}+\Cal E_\mu \subset \Cal E_\mu, 
$$
where we have assumed that the assertion is already proved with $\mu $
replaced by $\mu +1$. The spaces $\Cal E_\mu $ are thus determined by
local properties. Inside the set, the condition $u\in \Cal E_\mu $ only
means that $u$ is a $C^\infty $ function.

To determine the meaning of the condition $u\in \Cal E_\mu $ at a
boundary point we consider the case when $u$ has compact support in a
coordinate patch where $\Omega $ is defined by the condition $x_n>0$.

\example{Remark 2.2}
It will be useful to recall some formulas for power functions in one variable $t$ and
their Fourier transforms.
Denote as in \cite{H65} 
$$
I^\mu(t)
= \cases t^\mu /\Gamma (\mu +1)\text{ for }t>0,\\ 0 \text{ for }t\le 0,\endcases
\tag2.2$$ 
when $\operatorname{Re}\mu >-1$; it is called $\chi _+^\mu (t)$ in
\cite{H83}, Section 3.2. It is shown there  
that the distribution $I^\mu $
extends analytically from $\operatorname{Re}\mu >-1$ to $\mu \in{\Bbb C}$.
(For negative
integers, $I^{-k}= \delta _0^{k-1}$.) Moreover, \cite{H65} uses the
notation 
$(z ^\pm)^a$ for
the boundary values of $z^a$ from the half-planes
$\C_\pm=\{z\in{\Bbb C}\mid \operatorname{Im}z\gtrless 0\}$,
defined to be
real and positive on the positive real axis (they are denoted $(z\pm
i0)^{a}$ in \cite{H83}). Explicitly,
$$
(z^+)^a=\cases z^a\text{ for }z>0,\\ |z|^ae^{i\pi a}\text{ for
}z<0;\endcases
\quad
(z^-)^a=\cases z^a\text{ for }z>0,\\ |z|^ae^{-i\pi a}\text{ for
}z<0.\endcases \tag2.3
$$
Then, cf.\ \cite{H83}, Ex.\ 7.1.17, $I^\mu (t)$ 
has the Fourier
transform  
$$\F _{t\to\tau }I^\mu =e^{-i\pi (\mu +1)/2}(\tau ^-)^{-\mu -1}.\tag2.4$$  
We also note that when $\sigma >0$, translation by $-i\sigma
$ gives
$$
e^{-i\pi (\mu +1)/2}\F^{-1}(\tau -i\sigma )^{-\mu -1}=
\F^{-1}(\sigma +i\tau )^{-\mu -1}=
I^\mu e^{-t\sigma }.\tag2.5
$$
\endexample

\proclaim{Lemma 2.3}
An element $u\in \Cal E'({\Bbb R}^n)$ belongs to $\Cal E_\mu (\crnp)$,
if and
only if $u$ vanishes when $x_n<0$ and one can find $u_0,u_1,\dots \in
C_0^\infty ({\Bbb R}^{n-1})$ such that for every $N$
$$
\hat u(\xi )-\sum_0^{N-1}(\xi _n-i )^{-\mu -j-1}\hat u_j(\xi ')=O(|\xi
|^{-\operatorname{Re}\mu -N-1}), \; \xi \to\infty .\tag2.6
$$
Conversely, given such $u_0,u_1,\dots$ one can find $u\in\Cal E_\mu
(\crnp)$ satisfying this condition.

Here the argument of $\xi _n-i$ is chosen so that it tends to $0$ when
$\xi _n\to +\infty $.
\endproclaim

\demo{Proof} Any element $u\in \Cal E_\mu $ can be written
$u=v+\partial w/\partial x_n$ where $v$ and $w$ belong to $\Cal E_{\mu
+1}$. If the necessity of (2.6) has been proved when $\mu $ is
replaced by $\mu +1$ it follows therefore for $\mu $. Hence we may
assume that $\operatorname{Re}\mu >0$, thus $u=v x_n^\mu $ when
$x_n>0$, where $v\in C_0^\infty ({\Bbb R}^n)$. By forming a Taylor
expansion of $ve^{x_n}$ we can write for every $N$ 
$$
v=e^{-x_n}\sum_0^N v_j(x')x_n^j+R_N(x)
$$
where $v_j\in C_0^\infty ({\Bbb R}^{n-1})$ and $R_N(x)=O(x_n^N)$ when
$x_n\to 0$, $R_N(x)=O(e^{-x_n/2})$ when $x_n\to \infty $. Set
$R^0_N(x)=e^+r^+R_N(x)$. Then
$R^0_N(x)x_n^\mu $  has integrable derivatives of order $N$, so the
Fourier transform is $O(|\xi |^{-N})$. Now
$$
\hat u=\sum_0^\infty \hat v_j(\xi ')\F_{x_n\to\xi _n}(e^+r^+ e^{-x_n}x_n^{\mu +j })+
\F_{x\to\xi }( R^0_N(x)x_n^\mu ).
$$
By (2.5), $\F_{x_n\to \xi _n}(e^+r^+ e^{-x_n}x_n^{\mu +j })=\Gamma (\mu +j +1)e^{-i\pi (\mu +j +1)/2}(\xi
_n-i)^{-\mu -j-1}$,
so if we set 
$$
u_j=v_j \Gamma (\mu +j +1)e^{-\pi  i (\mu +j +1)/2},\tag2.7
$$
 it follows that (2.6) holds with the error term $O(|\xi
 |^{-N})$. Taking a few additional terms in the left hand side of
 (2.6) and noting that they can all be estimated in terms of the
 quantity on the right, we thus conclude that (2.6) is valid. 

On the other hand, if $u$ satisfies (2.6) we obtain with $v_j$
 defined by (2.7) that $u-e^{-x_n}\sum_0^{N-1}v_jx_n^{j+\mu }$ will
 be arbitrarily smooth if $N$ is large. This proves the sufficiency of
 (2.6). To prove the last statement we again assume that
 $\operatorname{Re}\mu >0$, take $\chi \in C_0^\infty ({\Bbb R})$ 
 equal to 1 when $|x_n|<1$ and define 
$$
u(x)=0,\, x_n\le 0,\quad u(x)=\sum_0^\infty e^{-x_n}v_j(x')x_n^{\mu
+j}\chi (x_na_j),\; x_n>0,
$$
where $a_j$ is chosen so large that the derivatives of the $j$th term
of order $\le j$ are all $\le 2^{-j}$. This is possible since
$(x_na_j)^\nu \chi ^{(k)}(x_na_j)$ is bounded uniformly in $x_n$ and
$a_j$ if $\operatorname{Re}\nu \ge 0$. This completes the proof.
\qed\enddemo

The particular case where $\mu $ is an integer is of special
importance. When $\mu \ge 0$ the space $\Cal E_\mu $ then consists of
all functions in $C^\infty (\comega)$ which vanish to the order $\mu $
at the boundary (that is, the derivatives of order $<\mu $ vanish
there), extrapolated by 0 outside. When $\mu <0$ we have the sum of a
function in  $C^\infty (\comega)$
  extrapolated as 0 in the complement of $\comega$, and multiple
 layers with $C^\infty $ densities and of order $<-\mu $ on
 $\partial\Omega $. This is the only case when $\Cal E_\mu $ contains
 elements supported  by $\partial\Omega $; in other words, the
 restriction of an element in $\Cal E_\mu $ to $\Omega $ determines it
 uniquely except when $\mu $ is a negative integer.

\example{Remark 2.4}
It was convenient in the proof of Lemma 2.1 to work with powers of
$\xi _n-i$ instead of powers of $\xi _n$, and one could also work
with powers of $\xi _n-i\sigma $ with a $\sigma >0$, e.g.\ $\sigma
=[\xi ']$; however $(\xi _n^-)^a$ are more 
convenient in some applications. In terms of these functions
we can rewrite (2.6) in
the form
$$
\hat u(\xi )-\sum_0^{N-1}(\xi _n^-)^{-\mu -j-1}\hat u'_j(\xi ')=O(|\xi
|^{-\operatorname{Re}\mu -N-1}), \; \xi \to\infty ,\; |\xi _n|>1,\tag2.8
$$
where $u'_j$ is a linear combination of $u_0,\dots,u_j$ with
coefficient 1 for $u_j$. Namely, insert  Taylor expansions
$(z-i)^a=(z^-)^a+(-i)a(z^-)^{a-1}+(-i)^2\frac12 a(a-1)(z^-)^{a-2}+\dots$ of the terms
$(\xi _n-i)^{-\mu -j-1}$, and regroup the resulting sums. 
Thus the $u'_j$ occurring in (2.8) are
in one to one correspondence with the $u_j$ in (2.6) and can be
chosen arbitrarily.

In particular, when $\mu =0$, so that $\E_\mu (\comega)=e_\Omega C^\infty (\comega )$,
$$
u_0=u'_0=-i\gamma _0u,\tag2.9
$$
where $\gamma _0u$ is the boundary value from $\Omega $.
\endexample

Consider a classical pseudodifferential operator $P$ in $\Omega_1$ of order $m \in
{\Bbb C}$.
Recall
 the notation for derivatives of the symbol in local coordinates: 
$$
p^{(\alpha )}_{(\beta
)}(x,\xi )=\partial_\xi ^\alpha \partial_x^\beta p(x,\xi ).\tag2.10
$$

The first question to investigate is when  $P$  maps $\Cal E_\mu $ into
$C^\infty (\comega )$ (more precisely, the restrictions to $\Omega $
belong to $ C^\infty (\comega )$). By the pseudo-local property of
$\psi $do's we
know that $Pu\in C^\infty (\Omega )$ for all $u\in \Cal E_\mu $. We
shall therefore only expect a restriction on $P$ at points on
$\partial\Omega $. Of course it is no restriction to assume $P$
compactly supported when studying a regularity problem.

\proclaim{Definition 2.5} A classical pseudodifferential operator of
order $m$ in $\Omega_1$ is said to satisfy the $\mu $-transmission condition
relative to $\Omega $ (in short: be of type $\mu $), when the symbol in any local coordinate
system 
satisfies
$$
{p_j}^{(\alpha )}_{(\beta )}(x,-N)=e^{\pi i(m-2\mu -j-|\alpha |)
}{p_j}^{(\alpha )}_{(\beta )}(x,N),\; x\in\partial\Omega ,\tag2.11
$$
for all $j,\alpha ,\beta $, where 
$N$ denotes the
interior normal of $\partial\Omega $ at $x$.
\endproclaim

\proclaim{Theorem 2.6} Let $P$ be a classical compactly supported
pseudodifferential operator of order $m$ in $\Omega_1$. In order that $r_\Omega Pu\in 
C^\infty (\comega )$ for all $u\in \Cal E_\mu (\comega)$, it is
necessary and sufficient that $P$ satisfies the $\mu $-transmission
condition.
\endproclaim

Since every polynomial satisfies this hypothesis with $\mu =0$ it
follows from the rules for coordinate changes that (2.11) is invariant under any change
of variables. In the proof of the theorem we may therefore use local coordinates such
that $\Omega $ is defined by the inequality $x_n>0$. The statement is
local, so it is enough to consider $Pu$ for $u\in \Cal E_\mu (\crnp)$
with compact support in the coordinate patch $U\subset {\Bbb
R}^n$. After modifying $P$ by an operator with symbol 0 we may assume
that $P$ is a compactly supported operator in $U$.

A key observation is the following elementary lemma.

\proclaim{Lemma 2.7} Let $q$ be a positively homogeneous function on
${\Bbb R}$ of degree $\sigma $, $\operatorname{Re}\sigma <-1$. For
$t>0$ we set $\varphi _\sigma (t)=t^{-\sigma -1}$ if $\sigma $ is not
an integer and $\varphi _\sigma (t)=t^{-\sigma -1}\log t$ if $\sigma $
is an integer. Then
$$
\int_{|\tau |>1}e^{it\tau }q(\tau )\, d\tau , \; t>0,
$$
is on $\rp$ equal to the sum of a function in $C^\infty
(\crp)$ and $C\varphi _\sigma (t)$. Here $C=0$ if and only if
$q(-1)=e^{i\pi \sigma }q(1)$, that is, if $q(\tau )= q(1)(\tau
^+)^\sigma $.
\endproclaim

\demo{Proof} Let $\gamma _+$ ($\gamma _-$) consist of
the real axis with the interval $(-1,1)$ replaced by a semi-circle in
the upper (lower) half plane. Then the two functions
$$
\int_{|\tau |>1}(\tau ^\pm)^\sigma  e^{it\tau }\, d\tau -\int_{\gamma ^\pm}(\tau ^\pm)^\sigma  e^{it\tau }\, d\tau 
$$
are integrals of $e^{it\tau } $ over semi-circles, hence obviously
entire analytic functions of $t$. By Cauchy's integral formula one
concludes that the integral over  $\gamma _+$ ($\gamma _-$) vanishes
for $t>0$ ($t<0$), and that it is homogeneous of degree $-\sigma -1$
 when $t<0$ ($t>0$). When $\sigma $ is not an integer, the two
functions $(\tau ^+)^\sigma $ and $(\tau ^-)^\sigma $ are linearly
independent, hence form a basis for positively homogeneous functions of
degree $\sigma $. This proves the lemma for non-integral $\sigma $.

To complete the proof it only remains to study
$$
\int_{|\tau |>1}(\tau ^\pm)^{\sigma -1}|\tau |\, e^{it\tau }\, d\tau 
$$
when $\sigma $ is an integer $\le -2$. When $\sigma =-2$ the last
integral is equal to 
$$
2\int_1^\infty \tau ^{-2}\sin t\tau \, d\tau =2t\int_{1/t}^\infty \tau
^{-2}\sin\tau \, d\tau .
$$
A Taylor expansion of $\sin\tau $ shows that the integral is equal to
$\log 1/t$ plus a function in $ C^\infty (\crp)$. This proves
the statement when $\sigma =-2$, and by successive integration it
follows for all integers $\sigma <-2$.
\qed\enddemo

\demo{Proof of Theorem {\rm 2.6}} 
Suppose that the theorem were already
proved with $\mu $ replaced by $\mu +1$. The necessity of (2.11) is
then obvious for it holds with $\mu $ replaced by $\mu +1$ and
$e^{-2\pi i}=1$. To prove its sufficiency we have to show that $PDu\in
 C^\infty (\comega )$ if $u\in \Cal E_{\mu +1}$ and $D$ is a
first order differential operator. Since $PDu=DPu+[P,D]u$ and $[P,D]$
satisfies (2.11) if $P$ does, the assertion follows. Hence we may
assume in what follows that $\operatorname{Re}\mu
>\operatorname{Re}m$. Then the product of $p(x,\xi )$ by the Fourier
transform of any compactly supported $u\in \Cal E_\mu (\crnp)$ is
integrable, so by an obvious regularization we obtain
$$
p(x,D)u=(2\pi )^{-n}\int p(x,\xi )\hat u(\xi )e^{ix\cdot\xi }\, d\xi .\tag2.12
$$
We shall introduce a Taylor expansion of $p$ in (2.12),
$$
p(x,\xi )=\sum_{|\alpha |<\nu }(\partial^{|\alpha |}p(x',0,0,\xi
_n)/\partial\xi ^{\alpha '}\partial x_n^{\alpha _n}) x_n^{\alpha
_n}\xi ^{\alpha '}/\alpha ! +\sum_{|\alpha |=\nu }r^\alpha (x,\xi )x_n^{\alpha _n}\xi ^{\alpha '},\tag2.13
$$
where 
$$
r^{\alpha }(x,\xi )=|\alpha |/\alpha !\int_0^1(1-t)^{|\alpha |-1}p^{(\alpha
')}_{(\alpha _n)}(x',tx_n,t\xi ',\xi _n)\, dt,
$$
where somewhat incorrectly we have used the notation $\alpha '$ for
$(\alpha ',0)$ and $\alpha _n$ for $(0,\alpha _n)$.
When $|\alpha '|> \operatorname{Re}m$ we can estimate $r^\alpha $ by
$(1+|\xi _n|)^{\operatorname{Re}m-|\alpha '|}$, and when $|\alpha
'|\le \operatorname{Re}m$ we can estimate by $(1+|\xi
|)^{\operatorname{Re}m-|\alpha '|}$ instead. Now we have
$$
\int r^\alpha (x,\xi )x_n^{\alpha _n}\xi ^{\alpha '}\hat u(\xi )e^{ix\cdot\xi
}\, d\xi =\int (i\partial_{ \xi _n})^{\alpha _n}\bigl(r^\alpha (x,\xi )\xi ^{\alpha '}\hat u(\xi )\bigr)e^{ix\cdot\xi
}\, d\xi .
$$
 Here the factor $x_n^{\alpha _n}$  was removed by an
integration by parts with respect to $\xi _n$ (using that $x_n^{\alpha
_n}e^{ix_n\xi _n}=(-i\partial_{\xi _n})^{\alpha _n}e^{ix_n\xi _n}$).
In view of (2.6) we conclude that the integral and its derivatives of
order $\le k$ are absolutely convergent, thus the integral defines a
$C^l$ function, provided that 
$$
l+\operatorname{Re}m-|\alpha '|-\alpha _n-\operatorname{Re}\mu <0.
$$
If we choose $\nu >k+\operatorname{Re}(m-\mu )$, the error term in
(2.13) will therefore only contribute a $C^l$ term to $p(x,D)u$.
The remaining problem is only to study the regularity of the partial
sums of the series obtained by replacing $p(x,\xi )$ by its Taylor
expansion in (2.12). Since $\hat u$ is rapidly decreasing when $\xi
\to\infty $ with $|\xi _n|<1$, this part of the integral in (2.12) is
infinitely differentiable. In view of (2.8) --- where we drop the
prime on $u'_j$ --- it only remains to examine when the partial sums
of the series
$$
\sum_{\alpha ,j,k }(2\pi )^{-n}\int_{|\xi _n|>1}{p_j}^{(\alpha
')}_{(\alpha _n)}(x',0,0,\xi _n)x_n^{\alpha _n}\xi ^{\alpha '}\hat
u_k(\xi ')(\xi _n^-)^{-\mu -k-1}e^{ix\cdot\xi }\, d\xi /\alpha !
$$
become arbitrarily smooth when the order of the sum goes to infinity. 
Here we can remove $x_n^{\alpha _n}$  by an
integration by parts with respect to $\xi _n$ as above.
The boundary terms
which then occur will give rise to only $C^\infty $ terms. Thus we are
reduced to examining the differentiability of the partial sums of the
series
$$
\sum_{\alpha ,j,k }D^{\alpha '}u_k(x')(2\pi )^{-1}\int_{|\xi
_n|>1}(i\partial_{\xi _n})^{\alpha _n}\bigl({p_j}^{(\alpha 
')}_{(\alpha _n)}(x',0,0,\xi _n)(\xi _n^-)^{-\mu
-k-1}\bigr)e^{ix_n\xi _n}\, d\xi _n /\alpha !.
$$
Since the functions $D^{\alpha '}u_k$ can be chosen arbitrarily in the
neighborhood of any point, or rather, linear combinations of them are
arbitrary, we conclude that for $P$ to have the required property it
is necessary and sufficient that for any $\alpha '$ and $k=0,1,\dots$
the partial sums of higher order of the series
$$
\sum_{\alpha _n,j }(2\pi )^{-1}\int_{|\xi _n|>1}(i\partial_{ \xi _n})^{\alpha _n}\bigl({p_j}^{(\alpha
')}_{(\alpha _n)}(x',0,0,\xi _n)(\xi _n^-)^{-\mu -k-1}\bigr)e^{ix_n\xi _n}\, d\xi _n/\alpha !
\tag2.14$$
are in $ C^\nu (\crp)=r^+C^\nu (\R)$ for any given $\nu $.
Here $(i\partial_{ \xi _n})^{\alpha _n}\bigl({p_j}^{(\alpha
')}_{(\alpha _n)}(x',0,0,\xi _n)(\xi _n^-)^{-\mu -k-1}\bigr)$ is
homogeneous of degree $m-j-|\alpha |-\mu -k-1$, so if $m-j-|\alpha
|-\mu -1=\sigma $, the degree is $\sigma -k$.
 
Now we shall apply Lemma 2.7. 
Noting that a finite sum $\sum c_j\varphi _{\sigma _j}(t) $ with
different $\sigma _j$ is in $ C^\nu (\crp)$ if and only if
$c_j=0$ when $-\sigma _j-1\le \nu $, we conclude that (2.14) has the
desired differentiability properties if and only if for each complex
number $\sigma $, each $\alpha '$ and $k=0,1,\dots$, each $x'$, the sum
$$
q(\xi _n)\equiv \sum_{m-j-|\alpha |-\mu -1=\sigma  }(i\partial_{ \xi _n})^{\alpha _n}\bigl({p_j}^{(\alpha
')}_{(\alpha _n)}(x',0,0,\xi _n)(\xi _n^-)^{-\mu -k-1}\bigr) /\alpha _n!\tag2.15
$$
is proportional to $(\xi _n^+)^{\sigma -k}$. (The sum of course
contains only finitely many terms.)

In view of the homogeneity of ${p_j}^{(\alpha ')}_{(\alpha
_n)}(x',0,0,\xi _n)$ of degree $m-j-|\alpha '|$, we
have for each term in the sum:
$$
\aligned
(i\partial_{ \xi _n})^{\alpha _n}&\bigl({p_j}^{(\alpha
')}_{(\alpha _n)}(x',0,0,\xi _n)(\xi _n^-)^{-\mu -k-1}\bigr) \text{
for $\xi _n>0$ equals}\\
&=i^{\alpha _n}\partial_{ \xi _n}^{\alpha _n}\bigl({p_j}^{(\alpha
')}_{(\alpha _n)}(x',0,0,1)\xi _n^{m-j-|\alpha '|-\mu -k-1}\bigr)
\\
&=i^{\alpha _n}{\ssize(m-j-|\alpha '|-\mu -k-1)\cdots(m-j-|\alpha |-\mu -k)}{p_j}^{(\alpha
')}_{(\alpha _n)}(x',0,0,1)\xi _n^{m-j-|\alpha |-\mu -k-1}
\\
&=i^{\alpha _n}{\ssize(m-j-|\alpha '|-\mu -k-1)\cdots(m-j-|\alpha |-\mu -k)}{p_j}^{(\alpha
')}_{(\alpha _n)}(x',0,0,1)\xi _n^{\sigma -k},
\endaligned\tag2.16
$$
whereas (cf.\ also (2.3))
$$
\aligned
(i\partial_{ \xi _n})^{\alpha _n}&\bigl({p_j}^{(\alpha
')}_{(\alpha _n)}(x',0,0,\xi _n)(\xi _n^-)^{-\mu -k-1}\bigr) \text{
for $\xi _n<0$ equals}\\
&=i^{\alpha _n}\partial_{ \xi _n}^{\alpha _n}\bigl({p_j}^{(\alpha
')}_{(\alpha _n)}(x',0,0,-1)|\xi _n|^{m-j-|\alpha '|-\mu -k-1}
e^{-\pi i(-\mu -k-1)}
\bigr)
\\
&=(-i)^{\alpha _n}{\ssize(m-j-|\alpha '|-\mu -k-1)\cdots(m-j-|\alpha |-\mu -k)}{p_j}^{(\alpha
')}_{(\alpha _n)}(x',0,0,-1)|\xi _n|^{\sigma-k }
e^{\pi i(\mu +k+1)}.
\endaligned$$
A function equal to (2.16) on $\R_+$ will be proportional to $(\xi
_n^+)^{\sigma -k}$ exactly when it on $\R_-$ has the value
$$
i^{\alpha _n}{\ssize(m-j-|\alpha '|-\mu -k-1)\cdots(m-j-|\alpha |-\mu -k)}{p_j}^{(\alpha
')}_{(\alpha _n)}(x',0,0,1)|\xi _n|^{\sigma -k}e^{\pi i(\sigma -k)}.
$$

Thus $q(\xi _n)$, where we for
fixed $\alpha '$, $k$, $\sigma $,  take the sum over 
$m-j-|\alpha
|-\mu -1=\sigma  $, 
is proportional to $(\xi _n^+)^{\sigma -k}$ if and
only if
$$
\multline
\sum_{m-j-|\alpha
|-\mu -1=\sigma   }
{\ssize(m-j-|\alpha '|-\mu -k -1)\dots(m-j-|\alpha |-\mu -k)}{p_j}^{(\alpha
')}_{(\alpha _n)}(x',0,0,1)e^{\pi i(\sigma  -k)}/ \alpha _n!=
\\
\sum{\ssize(m-j-|\alpha '|-\mu -k -1)\dots(m-j-|\alpha |-\mu -k)}(-1)^{\alpha _n}{p_j}^{(\alpha
')}_{(\alpha _n)}(x',0,0,-1)e^{\pi i(\mu +k+1)}/ \alpha _n!
.
\endmultline
$$
After the exponential factors have been moved to the same side and
integer powers of $e^{2\pi i}$ have been eliminated, we find
that $k$ occurs only in the polynomial factors, which are of degree
$\alpha _n$, all different. It follows that the coefficients have to
agree, that is  
$$
{p_j}^{(\alpha
')}_{(\alpha _n)}(x',0,0,1)e^{\pi i(m-j-|\alpha '|-2\mu )}={p_j}^{(\alpha
')}_{(\alpha _n)}(x',0,0,-1).\tag2.17
$$
This gives is a necessary and sufficient condition for $r^+P$ to map
$\Cal E_\mu(\crnp) $ into $ C^\infty (\crnp)$. But (2.17) is a consequence of (2.11),
and conversely, by differentiating (2.17) with respect to $x'$ and
using the homogeneity with respect to  $\xi _n$ we obtain
(2.11). This completes the proof of Theorem 2.6.
\qed\enddemo

Note that it suffices that the conditions in
(2.11) hold for the subset of derivatives ${p_{j}}^{(\alpha
')}_{(\alpha _n)}$ indicated in (2.17). A similar  sharpening is proved in \cite{GH90} for
more general,  not necessarily polyhomogeneous symbols,
in the case $\mu =0$.

In \cite{B69}, Boutet de Monvel with reference to the notes \cite{H65}
showed that (2.11) for $\psi $do's with analytic symbols implies a
mapping property as in Theorem 2.6 for functions analytic up to $\partial\Omega $. 

The product of two symbols of type $\mu _1$ resp.\ $\mu _2$ is clearly
of type $\mu _1+\mu _2$.

\example{Example 2.8} As simple examples, let us mention $(-\Delta
)^\nu $ and $\Lambda _\pm^\nu $ on $\R^n_+$ ($\nu \in{\Bbb C}$).
For $(-\Delta )^\nu $, of order $m=2\nu $,
the symbol $|\xi |^{2\nu }$ equals 1 for $\xi '=0$, $\xi _n=\pm 1$, so (2.11)
is satisfied with $\mu =\nu $; it is of type $\nu $.

For $\lambda _+^\nu   $, the principal symbol $(\lambda _{+}^\nu )_0$ is $(|\xi '|\overline\psi (\xi
_n/(a|\xi '|))+i\xi _n)^\nu $ (recall that $\psi (\pm\infty )=0$), so $(\lambda _{+}^\nu )_0(0,\pm
1)=(\pm i)^\nu $, satisfying (2.11) with $m=\nu $, $\mu =\nu $. The
difference between $\lambda _{\pm}^\mu $ and $(\lambda _{\pm}^\mu )_0$
is of order $-\infty $, since it has compact support in $\xi '$ and is rapidly
decreasing in $\xi _n$.
This shows that
 $\lambda ^\nu _+$ is
of type $\nu $. 

A similar study of $\lambda _-^\nu $ gives that it satisfies (2.11) with
$m=\nu $, $\mu =0 $, since the principal part clearly does so, and the remainder
is of order $-\infty $. Hence
 it is of type 0.

Moreover,  the modified symbols $\lambda _{\pm ,0}^{(\mu )}$, used in the
construction of order-reducing operators on a manifold (Theorem 1.3), are of type
$\mu $ resp.\ 0, since the exact symbols $\lambda ^\mu _\pm$ are used
near $\partial\Omega $, modulo smoothing terms.

\endexample

We also have, when $\Omega _1$ is compact:

\proclaim{Lemma 2.9} Let $A$ be a strongly elliptic second-order
differential operator with $C^\infty $-coefficients, and let $\nu \in{\Bbb C}$. Then 
the pseudodifferential operator $A^\nu  $ is of order $2\nu  $, and of
type $\nu $ for any smooth set $\Omega $.
\endproclaim

\demo{Proof} $A^\nu  $ is constructed by the method of Seeley
\cite{S67} (we recall that if 0 is an eigenvalue of $A$, $A^\nu $ is
taken zero
on the generalized eigenspace). First it is found that the resolvent $Q=(A-\lambda )^{-1}$
has the symbol in local coordinates
$$
\aligned
q(x,\xi ,\lambda )
&\sim \sum_{l\ge 0}q_{-l}(x,\xi,\lambda )
,\text{ where }
q_{0}=({a_0(x,\xi)-\lambda })^{-1},\\
 q_{-1}&={b_{1,1}(x,\xi
)}q_{0}^{2},\; \dots,\;
q_{-l}=\sum_{k=l/2}^{2l}{b_{l,k}(x,\xi)}{q_{0}^{k+1}},\; \dots\; ;
\endaligned
$$
with symbols $b_{l,k}$  independent of $\lambda $ and polynomial of
degree   $2k-l$ in $\xi$. (References are given e.g.\ in \cite{G96},
Remark 3.3.7.) The symbol of the $\nu $-th power of $A$ is
essentially constructed from this by a Cauchy integral together with $\lambda ^\nu $ around
the spectrum. The principal term gives $(a_0(x,\xi ))^\nu $, where, at
boundary points,
$$
a_0=s_0(x')\xi _n^2+O(|\xi _n||\xi '|)+O(|\xi '|^2),\quad s_0(x')\ne 0,
$$ with similar properties
as the Laplacian symbol above; the $\nu $-th power satisfies (2.11)
with
$m=2\nu $
and $\mu =\nu $. In the next
terms, when $q_0^{k+1}=c\partial_\lambda ^kq_0$ is inserted in the
integral and the $\lambda $-derivative is carried over to $\lambda
^\nu $, we get powers $(a_0(x,\xi ))^{\nu -k} $, that likewise satisfy
(2.11) with $\mu =\nu $, since the factors $a_0^{-k}$ are of type 0. It follows that $A^\nu $ is of type
$\nu $.
\qed
\enddemo

\example{Remark 2.10}
Consider $A$ as above and assume moreover that it has product structure near the boundary
$\partial\Omega $, i.e., coordinates can be chosen near $\partial\Omega $ such that
$A=D_n^2+A'(x',D')$ there with $A'$ strongly elliptic on $\partial\Omega $.
Then the associated Dirichlet-to-Neumann operator $P_{DN}$ (sending
$\gamma _0u$ to $\gamma _1u$ when $Au=0$) is
essentially a constant times $(A')^{\frac12}$, which is of order 1 and
type
$\frac12$ with respect to smooth subsets of $\partial\Omega $.
\endexample 

\example{Remark 2.11} 
When the equations (2.11) are satisfied with $\mu =0$ and
$m$ integer, they hold also if the normal vectors $N$ and $-N$
exchange roles. Then $P$ is of type 0 also for the exterior domain
$\Omega _1\setminus\Omega $; the so-called two-sided transmission
property.  This is the case treated in 
the Boutet de Monvel calculus.
\endexample

Noninteger transmission properties have been used in another context by
Hirschowitz and Piriou \cite{HP79} to investigate lacunas by
application of
Fourier integral operators; see also the survey by Boutet de Monvel \cite{B79}.

 \head 3. The Vishik-Eskin estimates \endhead

Consider a $C^\infty $ manifold $\Omega_1$, a relatively
compact subset $\Omega $ with $C^\infty $ boundary $\partial\Omega $,
and a classical pseudodifferential operator $P$ in $\Omega_1$. The operator
$P$ we assume to be {\it elliptic} in $\Omega_1$, that is, in a local
coordinate system where the symbol is $\sum p_j(x,\xi )$, the terms
being homogeneous of degree $m-j$, we have
$$
p_0(x,\xi )\ne 0\text{ for } 0\ne \xi \in{\Bbb R}^n.\tag3.1
$$
Further we assume that the $\mu $-transmission condition is fulfilled at least for $j=\alpha =\beta =0$, that is, we
assume that there is a number $\mu $ such that 
$$
p_0(x,-N) = e^{\pi i(m-2\mu )}p_0(x,N),\; x\in\partial\Omega , \tag3.2
$$
where $N$ denotes the interior normal of $\partial\Omega$ at $x$. If
$n>2$ the set $\{\xi \mid \xi \in{\Bbb R}^n, \xi \ne 0\}$ is simply
connected, so for fixed $x$ we can define $\log p(x,\xi )$ uniquely by
fixing the value at one point. When $n=2$, we impose this as a
condition on $p$, called the root condition in analogy with the
corresponding condition in the case of differential equations. Then we
have
$$
\log p_0(x,\xi +\tau N) - \log p_0(x,\tau N)= \log \big(p_0(x,\xi
+\tau N)/ p_0(x,\tau N)\big)\to0,\; \tau \to\infty .
$$
Hence 
$$
\log p_0(x,\xi +\tau N) - m\log|\xi |\to a_\pm(x),\; \tau \to\pm\infty , \tag3.3
$$
where $\exp a_\pm=p_0(x,\pm N)$. It follows from (3.2) that
$e^{a_-}=e^{\pi i(m-2\mu )+a_+}$, that is, $\mu \equiv
m/2+(a_+-a_-)/2\pi i$ (mod 1). We define the factorization index $\mu
_0$ by
$$
\mu _0=m/2+(a_+-a_-)/2\pi i,\tag3.4
$$
noting that for reasons of continuity this number, which is always
congruent to $\mu $, must be a constant on connected components of
$\partial\Omega $. 
(There is a remark in H\"ormander  \cite{H65} that much of the theory goes through
with light modifications
when $m$ and $\mu _0$ are allowed to be variable, referring to the
1964 Doklady notes preceding \cite{VE65, VE67}.)
Note that we may replace
$\mu $ by $\mu _0$ in (3.2).

We can now state the basic existence theorem for the Dirichlet
problem, due to Vishik and Eskin in the case $p=2$, cf.\ \cite{VE65, E81}, and extended to
$1<p<\infty $ by Shargorodsky \cite{S94}. 

\proclaim{Theorem 3.1} Let $P$ be elliptic of order $m$ satisfying
{\rm (3.2)} (and
the root condition if $n=2$), and assume the factorization index $\mu _0$ introduced
above to be constant on $\partial\Omega $. Then the mapping 
$$
\dot{H}^s_p(\comega )\ni u\mapsto r_\Omega  Pu\in\overline
H^{s-\operatorname{Re}m}_p(\Omega )\tag 3.5
$$
is a Fredholm operator if $s$ is a real number with
$1/p-1<s-\operatorname{Re}\mu _0<1/p$.
\endproclaim

In the proof one observes that it suffices to prove the a priori estimate
for smooth functions $$
\|u\|_{s}\le
C(\|r_\Omega Pu\|_{s-\operatorname{Re}m_0}+\|u\|_{s-1}),\; u\in
\dot{H}^s_p(\comega ),  \tag3.6
$$
together with an analogous estimate for the adjoint $\,^tP$. This
can be  reduced to the study of
``constant-coefficient'' symbols $p_0(x_0,\xi )$ for
$x_0\in\partial\Omega $ in the case $\Omega
=\rnp$. Here there is a factorization 
$$
p_0(x_0,\xi )=p_-(x_0,\xi
)p_+(x_0,\xi )\tag3.7
$$ 
with $p_\pm$ of degree $\mu _0$ resp.\ $m-\mu _0$, extending as
analytic functions of $\xi _n$
to ${\Bbb C}_-$ resp.\ ${\Bbb C}_+$, hence
defining operators preserving support in $\crnp$ resp.\
$\crnm$. Details on the factorization and its application to obtain
the estimates are found e.g.\ in \cite{E81} \S6, 7, 19, extended to
$L_p$-spaces in  \cite{S94}. (See
(1.10)ff.\ concerning sign conventions.) Those works moreover treat
systems $P$ and cases where $\mu _0$ depends on $x\in\partial\Omega $;
then the interval where $s$ runs has a smaller length.

\example{Example 3.2} 
When $A^\nu $ is defined as in Lemma 2.9, the principal symbol at a
boundary point $(x',0)$ has the factorization
$$
a_0(x',0,\xi ',\xi _n)^\nu =s_0(x')^\nu (m^+(x',\xi ')-\xi _n)^\nu (m^-(x',\xi ')-\xi _n)^\nu ,
$$
where $m^\pm$ are the roots in ${\Bbb C}_\pm$, respectively, of the
characteristic polynomial of degree 2. Here $(m^\pm(x',\xi ')-\xi
_n)^\nu $ extends analytically to ${\Bbb C}_\mp$, respectively. Thus
the factorization index equals $\nu $, and Theorem 3.1 applies with
$s-\operatorname{Re}\nu \in \,]-1/p',1/p[\,$.

Let $\nu =a\in\rp$. In the application of the theorem, $s\in
a+\,]-1/p',1/p[\,$, so regardless of how regular $r_\Omega Pu$ is, this
gives at best $u\in \dot H^{a+1/p-0}_p(\comega)$. 
When $p>n/a$, Sobolev embedding gives  $u\in
C^{a+1/p-n/p-0}(\comega)$ with boundary value zero. For
$p\to\infty $ we get $u\in C^{a-0}(\comega)$. It is pointed out in
Ros-Oton and Serra \cite{RS14} for $(-\Delta )^a$ with\ $a\in
\,]0,1[\,$ that 
the exponent 
$a-0$ cannot in general be lifted to values $>a$.

There are similar considerations for strongly elliptic $2m$-order
differential operators. Here the principal symbol at the boundary
factors into two polynomials in $\xi _n$ of degree $m$
with roots in $\C_\pm$, respectively. The $\nu $'th power is then of
order $2\nu m$ and type $\nu m$, and has factorization index $\nu m$.

More generally, let $P$ be of order $m\in{\Bbb C}$ with an {\it even}
symbol, that is, $p_j(x,-\xi )=(-1)^jp_j(x,\xi )$ for all $j\ge
0$. Then in view of the homogeneity of each $p_j$, $p$ satisfies
(2.11) with $\mu =m/2$. For the principal symbol, $a_+=a_-$ in
(3.3)--(3.4), so the factorization index is $m/2$. (One can also
include a skew factor $e^{i\pi \varrho }$.) 

The integral operators treated in the recent work of Ros-Oton and
Serra \cite{RS14} have these properties, with $m=2s$, when the kernel is smooth (outside
0).

Note that (2.11) is only required for the interior normal $N(x)$ to a
given smooth subset $\Omega $; the above examples have the property
with respect to {\it all directions}.

\endexample

The new task is to characterize the regularity of $u$ when
$P u$ is given in more smooth spaces. There is a preparatory result in
\cite{H65} on
``tangential regularity'' which follows by classical arguments due to
Nirenberg.

Let $\Omega $ be the half ball $\{ x\in{\Bbb R}^n\mid |x|<1,x_n>0\}$.
 The unit ball we denote by $\widetilde\Omega $. By 
$\dot{H}_{p,\operatorname{loc}}^{s}(\Omega ')$ and $ 
 \overline
H_{p,\operatorname{loc}}^{s-\operatorname{Re}m}(\Omega
 )$ we denote the distributions
which multiplied with functions in $C_0^\infty (\widetilde \Omega )$
give elements in the analogous spaces in $\rnp$. Here $\Omega '=\{
 x\in{\Bbb R}^n\mid |x|<1,x_n\ge 0\}$.

\proclaim{Theorem 3.3} Let $P$ satisfy the hypotheses of Theorem
{\rm 3.1}. If $-1/p'<s-\operatorname{Re}\mu <1/p$ and $t_0,t_1$ are
real numbers, then
$$
u\in\dot{H}_{p,\operatorname{loc}}^{s,t_0}(\Omega '), 
r_\Omega Pu\in \overline
H_{p,\operatorname{loc}}^{s-\operatorname{Re}m,t_1}(\Omega )\tag3.8
$$
implies that
$$
u\in\dot{H}_{p,\operatorname{loc}}^{s,t_1}(\Omega '),\tag3.9 
$$
\endproclaim

\demo{Proof} It is no restriction to assume that $t_1-t_0$ is a
positive integer, for we may always decrease $t_0$. It suffices to
prove the theorem when $t_1-t_0=1$. Now we claim that for every
compact subset $K$ of $\Omega $, and every real number $t$ there is a
constant $C$ such that 
$$
\|u\|_{s,t}\le C(\|r_\Omega Pu\|^{\;}_{s-\operatorname{Re}
m,t}+\|u\|_{s-1,t})\tag 3.10
$$
for all $u\in C_0^\infty (K)$, hence for all $u\in \dot H^{s,t}$ with
support in $K$. In fact, this follows from by applying (3.6) to $[D']^tu$, cut off conveniently. We may
replace the last term in (3.10) by the larger quantity
$\|u\|_{s,t-1}$. Now assume that (3.8) is fulfilled with $t_0=t$,
$t_1=t+1$. Then $\varphi u$ satisfies the same hypothesis if $\varphi
\in C_0^\infty (\widetilde \Omega )$. 
Let therefore $u$
have compact support in $\Omega '$. Denote by $u_h$ the convolution of
$u$ by the Dirac measure at $(h_1,\dots,h_{n-1}, 0)=h$, that is, $u_h$
is a tangential translation of $u$. Let $P_h$ be the analogous
translation of $P$. Then 
$$
P(u_h-u)/|h|=(f_h-f)/|h|+ (P-P_h)/|h|\; u, \tag3.11
$$
where $f=Pu$. Since
$$
\|(f-f_h)/|h| \|^{\;}_{s,t}\le \| f\|^{\;}_{s,t+1},
$$
 and since $(P-P_h)/|h|$ is continuous from $H^{s,t}_p$ to
 $H^{s-\operatorname{Re}m,t}_p$ uniformly when $h\to 0$, we conclude
 using (3.10) that $\|(u_h-u)/|h|\|_{(s,t)}$ is bounded when $h\to
 0$. Hence $\|D_ju\|_{(s,t)}<\infty $ when $j<n$, which proves that
 $u\in \dot H_{(s,t+1)}$. \qed
\enddemo

\head 4. Solvability of homogeneous problems \endhead

For the study of solvability, we first set the $H_p^{\mu (s)}$-spaces in relation to $\E_\mu $.
In the following
we assume that $\comega $ is compact, unless otherwise mentioned.

\proclaim{Proposition 4.1} $1^\circ$ Let $s>\operatorname{Re}\mu -1/p'$. For
any compact $K$, $u\in \Cal E_\mu (\crnp)\cap \Cal E'(K)$ implies
 $u\in H^{\mu (s)}_p(\crnp)$. Similarly, $
\E_\mu(\comega)
\subset H^{\mu (s)}_p(\comega)$.

$2^\circ$ We have that $\bigcap _{s}H^{\mu
(s)}_p(\crnp)\subset\E_\mu (\crnp)$, and that $$
\bigcap _{s}H^{\mu (s)}_p(\comega)=\E_\mu
(\comega).\tag 4.1$$ 

$3^\circ$ Moreover, $\E_\mu (\crnp)\cap \dot{\E}'(\crnp)$, resp.\ $\E_\mu (\comega)$,
is dense in $H_p^{\mu (s)}(\crnp)\cap \dot{\E}'(\crnp)$ resp.\ $H_p^{\mu
(s)}(\comega)$, when $s>\operatorname{Re}\mu -1/p'$.
\endproclaim

\demo{Proof} $1^\circ$. Let $u\in \Cal E_\mu (\crnp)\cap \Cal
E'(K)$. Then by (2.8), we
have for  $|\xi _n|>1$, $M\in{\Bbb N}$,  and any $N$, 
$$
\hat u(\xi )=\sum_{j=0}^{M -1} \hat u_j(\xi ')(\xi _n^-)^{-\mu -j-1}+O([\xi ']^{-N}|\xi _n|^{-\operatorname{Re}\mu-M -1}),\tag4.2
$$
where the $\hat u_j$ are in $\SD(\R^{n-1})$. To estimate $\Xi _+^\mu
u$, we shall calculate
$\hat u(\xi )(\xi _n-i[\xi '])^\mu $, where we note that $(\xi
_n-i[\xi '])^\mu 
=(-i)^\mu ([\xi ']+i\xi _n)^\mu =(-i)^\mu \chi _+^\mu 
$.
There are Taylor expansions (for $|\xi _n|>1$, say)
$$
\multline
(\xi _n-i[\xi '] )^{\mu }=(\xi _n^-)^\mu +c_1[\xi '] (\xi
_n^-)^{\mu -1}+\dots+c_{l-1}[\xi '] ^{l-1}(\xi _n^-)^{\mu -l+1}\\
+O([\xi '] ^{l+[\operatorname{Re}\mu -l]_+}|\xi _n|^{\operatorname{Re}\mu -l}).
\endmultline\tag4.3
$$
Insertion gives (with $c_0=1$):
$$
\aligned
&\F(\Xi _+^\mu u)(-i)^{\mu }=\hat u(\xi )(\xi _n-i[\xi '] )^\mu \\
&=\sum_{j=0}^{M-1}\hat u_j(\xi
')(\xi _n^-)^{-\mu -j-1}\Bigl[\sum _{l=0}^{M-j-1}c_l[\xi ']^l  (\xi
_n^-)^{\mu -l}+O([\xi '] ^{M-j+[\operatorname{Re}\mu -M+j]_+}|\xi _n|^{\operatorname{Re}\mu -M+j})\Bigr]\\
&\quad +O([\xi ']^{-N}|\xi _n|^{-M-1})\\
&=\sum_{j=0}^{M-1}\sum_{l=0}^{M-j-1}\hat u_j(\xi
')c_l[\xi '] ^l(\xi _n^-)^{-j-l-1}
+O([\xi ']^{-N}|\xi _n|^{-M-1})\\
&=\sum_{j=0}^{M -1} \sum_{k=0}^jc_{jk}\hat u_k(\xi ')[\xi ']
^{j-k}\xi _n^{ -j-1}+O([\xi ']^{-N}|\xi _n|^{-M-1}).
\endaligned\tag4.4
$$
In the last step we replaced $l,j$ by $j'=l+j$ and $k'=j$, and
removed the primes. 
The $c_{kj}$ are constants, with
$c_{jj}=1$. (It is also for later purposes that we account for
this in detail.)

The terms in the sum are
Fourier transforms of functions in $\ol \SD (\rnp)$, and the remainder
is bounded by $\ang \xi ^{-N'}$ for $N'\le \min\{N,M+1\}$, so by letting
$N,M\to\infty $, we see that any $\ol
H^t_p(\rnp)$-norm of $\Xi _+^\mu u$ is bounded. 

The result for $\comega$ follows by using the above in local
coordinate patches where $d(x)=x_n$.

$2^\circ$. Now let  $u\in\bigcap _{s}H^{\mu
(s)}_p(\crnp)$. Then $v=r^+\Xi _+^\mu u\in \bigcap _{t}\ol
H^{t}_p(\rnp)$, which consists of $C^\infty (\crnp)$-functions with all
$L_p$-norms of derivatives bounded. In view of Proposition 1.7,
$u=\Xi ^{-\mu }_+e^+v
$. By Lemma 2.1, $v$ has an expansion as in (2.8) with $\mu =0$,
and the multiplication by $([\xi ']+i\xi _n)^{-\mu }=i^{-\mu }(\xi
_n-i[\xi '])^{-\mu }$ gives a function with an expansion (2.8)
with the actual $\mu $, so
we conclude from Lemma 2.1 with (2.8) that $u\in \E_\mu (\crnp)$. 

For $\comega$ we find from this by localization that  
$\bigcap _{s}H^{\mu
(s)}_p(\comega)\subset\E_\mu (\comega)$; here there is equality in
view of $1^\circ$.

$3^\circ$. To show
that $\E_\mu \cap \dot{\E}'(\crnp)$ is dense in the set of all $u\in \dot{\SD}'(\crnp)$
satisfying (1.24),
we first take a sequence $v_j\in  C^\infty (\crnp)$ of compactly
supported functions approximating
$\Cal F^{-1}(\xi _n-i[\xi '])^\mu  \hat u$ in the norm
$\|\;\|_{\ol H_p^{s-\operatorname{Re}\mu }},$ and also in
the topology of $\Cal S$ outside a neighborhood of $\supp u$ (which is
possible since the function to approximate agrees with a function in
$\SD$ there). 
Define $v_j=0$ in $\rnm$. Set
$u_j=\Cal F^{-1}((\xi _n-i[\xi '] )^{-\mu } \hat v_j)$.  This is an element of $\E_\mu $ in view of Lemma
2.2 (the Fourier transform is the product of that of $v_j$ and
$(\xi _n-i[\xi '] )^{-\mu }$, and the behavior of the Fourier
transform of $v_j$ is described by Lemma 2.3 with $\mu =0$). Then by
Proposition 1.7,
$u_j\to u$ in the norm in (1.24), and also in the topology of $\SD$
outside a neighborhood of $\supp u$. Hence we can cut off  $u_j$ there
without disturbing the convergence in order to obtain an approximating
sequence with compact supports. 

The statement for $H_p^{\mu (s)}(\comega)$ follows by localization.
\qed  
\enddemo

In the next theorems we use the order-reduction operators to reach
situations where we can draw on results from the Boutet de Monvel
calculus. The calculus was established in \cite{B71} and is moreover presented
in detail e.g.\ in \cite{G96, G09}, see also \cite{G90}.

\proclaim{Theorem 4.2} Let the $\psi $do $P$ on $\R^n$ be of
order $m$, and type $\mu $ relative to $\rnp$, and compactly supported. Then for
$s>\operatorname{Re}\mu -1/p'$ and $u\in H^{\mu (s)}_p(\crnp)$,
$$
\|r^+Pu\|^{\;}_{\ol H_p^{s-\operatorname{Re}m}}\le C \|r^+\Xi _+^\mu
u\|^{\;}_{\ol H_p^{s-\operatorname{Re}\mu }}, \quad
\|r^+Pu\|^{\;}_{\ol H_p^{s-\operatorname{Re}m}}\le C' \|r^+\Lambda _+^\mu
u\|^{\;}_{\ol H_p^{s-\operatorname{Re}\mu }}.
\tag4.5
$$

Similarly, for a $\psi $do $P$ on the manifold $\Omega _1$ of order
$m$, and type $\mu $ on $\Omega$, one has for $u\in H^{\mu (s)}_p(\comega)$,
$$
\|r_\Omega Pu\|^{\;}_{\ol H_p^{s-\operatorname{Re}m}(\Omega )}\le C
\|r_{\Omega }\Lambda _+^{(\mu )} u\|^{\;}_{\ol H_p^{s-\operatorname{Re}\mu }(\Omega )}.\tag4.6
$$

In other words, $r^+P$ maps $H_p^{\mu (s)}$ continuously into $\ol
H_p^{s-\operatorname{Re}m}$ when $s>\operatorname{Re}\mu -1/p'$.
\endproclaim

\demo{Proof}
By definition,
$v=r^+\Lambda _+^\mu u\in \ol
H^{s-\operatorname{Re}\mu }$, and by Proposition 1.7, 
$u=\Lambda _+^{-\mu }e^+v$ then. Thus we can write 
$$
r^+Pu=r^+P\Lambda _+^{-\mu }v.
$$
Moreover, by (1.15),
$$ 
\|r^+Pu\|^{\;}_{\ol H_p^{s-\operatorname{Re}m}}\simeq \|\Lambda ^{\mu
-m}_{-,+}r^+Pu\|^{\;}_{\ol H_p^{s-\operatorname{Re}\mu }},
$$
where 
$$
\Lambda ^{\mu -m}_{-,+}r^+Pu=r_+\Lambda _-^{\mu -m}Pu
$$
in view of Remark 1.1 (since the action of $\Lambda ^{\mu -m}_{-,+}$
is independent of how $r_+Pu$ is extended). Altogether,
$$
\|r^+Pu\|^{\;}_{\ol H_p^{s-\operatorname{Re}m}}\simeq
\|r^+Qe^+v\|^{\;}_{\ol H_p^{s-\operatorname{Re}\mu }},\text{ where }Q=\Lambda _-^{\mu -m}P\Lambda _+^{-\mu }.
$$
Here $Q$
is of order 0 and type
0, hence belongs to the Boutet de Monvel calculus (as noted in
Remark 2.11), and we have from \cite{G90} that $Q_+=r^+Qe^+$ is
continuous from $\ol H_p^{s-\operatorname{Re}\mu }$ to itself, since
$s>\operatorname{Re}\mu -1/p'$. This implies the second inequality in
(4.5), and the first one follows in view of  Lemma 1.2.

For $\Omega $ we obtain the result either by using the above in local
coordinates or by repeating the proof using $\Lambda _\pm^{(\mu )}$.
\qed
\enddemo

In  the notes \cite{H65}, the proof of this theorem for $p=2$ takes up
much space 
and involves a number of other tricks, needed because
 the order-reducing operators $\Lambda ^\mu _\pm$ were not known
then. Finiteness of {\it all seminorms} 
$$
u\mapsto \|r_\Omega Pu\|^{\;}_{\ol H_2^{s- \operatorname{Re}m}(\Omega )}, \tag4.7
$$
with $P$ of type $\mu $ and any order $m$, was taken as the {\it definition} of
the topology of $H_2^{\mu (s)}(\comega)$, and a large effort went into
showing that on $\rnp$, finiteness of $\|r^+\Xi ^\mu
_+u\|^{\;}_{\ol H_2^{s-\operatorname{Re}\mu }}$ suffices, or more precisely, finiteness of $\|r^+(\ang{D'}+\partial_n)^\mu
u\|^{\;}_{\ol H_2^{s-\operatorname{Re}\mu }}$ suffices. It comes in as a
special case when (4.7) is investigated for $P=(1-\Delta )^\mu $, $m=2\mu $.

The mapping property was proved for operators of type 0 and
any real order $m$ in \cite{GH90} for $L_2$-spaces, including more
general, not polyhomogeneous symbols in $S^m_{\varrho ,\delta }$.
 (This covers classical symbols of order $m\in{\Bbb C}$ and type 0, since they
 are in $S^{\operatorname{Re}m}_{1,0}$.)

\proclaim{Proposition 4.3} Let $s>\operatorname{Re}\mu-1/p'$. Both for
spaces over $\crnp$ and over $\comega$, we have that 
$$
H^{\mu (s)}_p\subset H_p^{(\mu -1)(s)},\tag4.8
$$
and the norms are equivalent on $H_p^{\mu (s)}$.

\endproclaim

\demo{Proof}
When $u\in H^{\mu (s)}_p(\crnp)$ for some $s>\operatorname{Re}\mu
-1/p'$, then
$$
\aligned
\|r^+\Xi ^{\mu -1}_+u\|^{\;}_{\ol H_p^{s-\operatorname{Re}\mu +1}}&\simeq
\sum_{j\le n}\|D_jr^+\Xi ^{\mu
-1}_+u\|^{\;}_{\ol H_p^{s-\operatorname{Re}\mu }}\\
&=
\sum_{j\le n}\|r^+D_j\Xi ^{\mu
-1}_+u\|^{\;}_{\ol H_p^{s-\operatorname{Re}\mu }}\le C \|r^+\Xi ^{\mu
}_+u\|^{\;}_{\ol H_p^{s-\operatorname{Re}\mu }},
\endaligned
$$
where we could use Theorem 4.2 in the last step, since $D_j\Xi _+^{\mu
-1}$ is of type $\mu $ and order $\mu $.
On the other hand, since $r^+\Xi _+^\mu u=r^+([D']+iD_n)\Xi _+^{\mu
-1}u$,
$$
\|r^+\Xi _+^\mu u\|^{\;}_{\ol H_p^{s-\operatorname{Re}\mu }}\le C\|r^+\Xi
_+^{\mu -1}u\|^{\;}_{\ol H_p^{s-\operatorname{Re}\mu +1}}. 
$$
Altogether, (4.8) holds, with equivalent norms on  $H_p^{\mu
(s)}$. Moreover, ${\Xi }^\mu _+$ can be replaced by $\Lambda _+^\mu $
in the inequalities in view of Lemma 1.2.

The
statements carry over to the manifold situation by localization.
\qed\enddemo

The $H^{\mu (s)}$-spaces serve the purpose of describing the
regularity of solutions with data in more regular Sobolev spaces than
the result of Vishik and Eskin (Theorem 3.1) allows.
We can now show the main regularity result for homogeneous boundary
problems (proved for $p=2$ in \cite{H65}), obtaining moreover a
formula for a parametrix:

\proclaim{Theorem 4.4} Let $P$ be classical elliptic of order $m\in{\Bbb C}$
on $\Omega _1$ and of type $\mu _0\in{\Bbb C}$ relative to $\Omega $, and  
with factorization index $\mu _0$. Let
$s>\operatorname{Re}\mu _0-1/p'$, and let $u\in \dot H_p^\sigma (\comega) $ for
some $\sigma >\operatorname{Re}\mu _0-1/p'$. If $r^+ Pu\in \ol H_p^{s-\operatorname{Re}m_0}(\Omega )$,
then $u\in H_p^{\mu _0(s)}(\comega)$. The mapping 
$$
H^{\mu _0(s)}_p(\comega )\ni u\mapsto r^+  Pu\in\overline
H^{s-\operatorname{Re}m}_p(\Omega )\tag 4.9
$$
is Fredholm, and has the parametrix 
$$
R=\Lambda ^{(-\mu
 _0)}_+e^+  \widetilde {Q_+}\Lambda ^{(\mu _0-m)}_{-,+}\colon 
\ol H^{s-\operatorname{Re}m}_p(\Omega  )
\to H^{\mu _0(s)}_p(\comega ) ,  \tag4.10
$$ 
where $\widetilde { Q_+}$ is a parametrix of $Q_+=r^+Qe^+$, with
$$
Q=\Lambda ^{(\mu _0-m)}_{-} P\Lambda ^{(-\mu _0)}_+,\tag4.11$$ 
elliptic of order
and type $0$, with factorization index $0$.

In particular, if $r^+ Pu\in C^\infty (\comega)$, then  $u\in \E_{\mu
_0}(\comega)$. The mapping 
$$
\E_{\mu _0}(\comega)\ni u\mapsto r^+ Pu\in C^\infty (\comega)\tag4.12
$$
is Fredholm.
\endproclaim

\demo{Proof} 
Note first that there is a $\sigma _0\le \min \{s,\sigma \}$ with $\sigma _0\in
\operatorname{Re}\mu _0+\,]-1/p',1/p[\,$. Theorem 3.1 (by
Vishik-Eskin-Shargorodsky) applies with
$s$ replaced by $\sigma _0$ to show the Fredholm solvability of
$r^+Pu=f\in \ol H_p^{\sigma _0-\operatorname{Re}m}$ with solution $u\in
\dot H^{\sigma _0}_p$. We must show that this solution lies in
$H^{\mu _0(s)}$. It already lies in $H^{\mu _0(\sigma _0)}$, since
$\Lambda _+^{(\mu _0)}u\in \ol H_p^{\sigma _0-\operatorname{Re}\mu
_0}\subset \dot H_p^{-1/p'+0}$.

To discuss the solvability of 
$$
r^+ Pu=f\in \ol H_p^{s-\operatorname{Re}m}(\Omega ),\tag 4.13
$$
in spaces with general $s$ we prefer to start from scratch, using devices from Theorem 4.2. Compose to the left with $\Lambda ^{(\mu _0-m)}_{-,+}$; this gives the equivalent
problem
$$
\Lambda ^{(\mu _0-m)}_{-,+}r^+ Pu=g,\text{ where }g=\Lambda ^{(\mu _0-m)}_{-,+}f\in \ol H_p^{s-\operatorname{Re}\mu _0}(\Omega ),\tag 4.14
$$
when we recall (1.20). Note that
$f=\Lambda ^{(m-\mu _0)}_{-,+}g$. Moreover, in view
of Remark 1.1,
$$
\Lambda ^{(\mu _0-m)}_{-,+}r^+ Pu=r^+\Lambda ^{(\mu _0-m)}_{-} Pu.
$$

Now set
 $v=r^+\Lambda ^{(\mu _0)}_+u$; then $u=\Lambda ^{(-\mu
 _0)}_+e^+v$ 
by Proposition 1.7. Expressed in terms of $g$ and $v$, equation (4.13) becomes 
$$
Q_+v=g; \quad g\text{ given in }\ol H_p^{s-\operatorname{Re}\mu _0}(\Omega ),\tag4.15
$$
where we have defined $Q$ by (4.11).
 
The properties of $P$ imply that $Q$ is elliptic of order 0 and type 0 and has
factorization index 0; in particular, it belongs to the Boutet de
Monvel calculus.
The principal symbol at the boundary $q(x',0,\xi )$ has a factorization
$q=q^+q^-$, in symbols $q^{\pm}(x',\xi )$ of plus/minus type and order
0. (We here use upper indices $\pm$ to avoid confusion with the lower plus-index
for truncation.) The associated operators on $L_2({\Bbb R})$ satisfy
$$
q_+=r^+q(x',0,\xi ',D_n)e^+=r^+q^-(x',\xi ',D_n)(e^+r^++e^-r^-)q^+(x',\xi ',D_n)e^+=q^-_+q^+_+,\tag4.16
$$
since $r^+q^-e^-$ and $r^-q^+e^+$ are zero. Let $\tilde q(x',\xi
)=1/q(x',\xi )$, it likewise has a factorization $\tilde q=\tilde
q^+\tilde q^-$ in plus/minus symbols, with $\tilde q^\pm=1/q^\pm$.
Now for the associated operators on ${\Bbb R}$,
$$
r^+q^+e^+r^+\tilde q^+e^+=r^+q^+\tilde q^+e^+-r^+q^+e^-r^-\tilde q^+e^+=I_{\rp},
$$
since $\tilde q^+$ preserves support in $\crp$ so that $r^-\tilde
q^+e^+=0$. One checks similarly that \linebreak $r^+\tilde q^+e^+r^+
q^+e^+=I_{\rp}$, and that also $r^+q^-e^+r^+\tilde q^-e^+=I_{\rp}$,
$r^+\tilde q^-e^+r^+ q^-e^+=I_{\rp}$. In other words,
$$
q^\pm(x',\xi ',D_n)_+\text{ has the inverse }\tilde q^\pm(x',\xi
',D_n)_+\text{ in }L_2(\rp).\tag4.17
$$
In view of (4.16), $q(x',0,\xi ',D_n)_+$ therefore has the inverse
$$
\widetilde {q_+}=(q^-_+q^+_+)^{-1}=\tilde q^+_+\tilde q^-_+.\tag4.18
$$
(More precisely, with notation from the Boutet de Monvel calculus as
described e.g.\ in \cite{G09} p.\ 284ff., $\widetilde {q_+}=\tilde
q_+-L(\tilde q^+,\tilde q^-)$, where the singular Green operator  $L(\tilde q^+,\tilde
q^-)=g^+(\tilde q^+)g^-(\tilde q^-)$ is generally nonzero.)

We see that $r^+q(x',0,\xi
',D_n)e^+$ is invertible as a boundary symbol operator,
and thus $Q_+=r^+Qe^+$ defines an elliptic boundary problem (without auxiliary trace or
Poisson operators) in the Boutet de Monvel calculus, hence defines a Fredholm
operator in $\dot H^t _2(\comega)=\ol H^t
_2(\Omega )$ for $|t |<\frac12$. (This is also shown in Vishik and Eskin's theorem, cf.\ Theorem 3.1.)

By [G90], $Q_+$ is continuous in  $\ol H_p^t(\Omega )$ for $t>-1/p'$.
We shall denote by $\widetilde {Q_+}$  a parametrix of $Q_+$;
likewise continuous in $\ol H_p^t(\Omega )$ for $t>-1/p'$. (It will
in general be of the form $\widetilde Q_++G$, where $\widetilde Q$ is
a parametrix of $Q$ and $G$ is a singular Green operator.) Thus,
 solutions of $Q_+v=g$ with in $g\in \ol H_p^t(\Omega )$
for some $t>-1/p'$ are in $\ol H_p^t(\Omega )$, and 
$$
Q_+\colon \ol H^t_p(\Omega )\to \ol H^t_p(\Omega )\text{ is Fredholm for all }t>-1/p'.
$$

It follows that the solutions of (4.15) satisfy $v\in \ol
H_p^{s-\operatorname{Re}\mu _0}(\Omega )$, so the solutions of the
original problem (4.13) satisfy $u\in H_p^{\mu _0(s)}(\Omega
)$. Retracing the steps, we find that (4.10) is a parametrix of $r^+P$.
The Fredholm property also follows.

Finally, the solvability with right-hand side in $C^\infty (\comega )$ is
deduced from the above by use of Proposition 4.1.\qed  
\enddemo

The proof in \cite{H65} of the Fredholm property in the $L_2$-case was
based on Theorems 3.1 and  3.3
together with certain intricate results on ``partial hypoellipticity at the
boundary'' (valid for general $P$ of type $\mu $ for which $\partial\Omega $ is
non-characteristic).

\example{Example 4.5}
Let us check how this looks in the well-known case of the
Laplace-Beltrami operator,
$P=\Delta $. It is of order 2 and type 0, and has factorization index 1
(cf.\ Example 3.2). Let $s>1-1/p'=1/p$, so $f$ is given in $\ol
H^{s-2}_p$ with $s-2>-2+1/p$. From Example 1.6 with $m=1$ we have that
$H_p^{1(s)}=\{u\in \ol
H^s_p\mid \gamma _0u=0\}$.
Thus $u$ is the solution of the homogeneous Dirichlet problem:
  $\Delta u=f$ in $\Omega $, $\gamma _0u=0$.
\endexample

\example{Remark 4.6}
Not all elliptic $\psi $do's $P$ of order and type 0 have $P_+$
elliptic without supplementing trace or Poisson operators. For
example, $P=\Lambda ^{(1)}_-\Lambda ^{(-1)}_+$ has $P_+=\Lambda
^{(1)}_{-,+}\Lambda ^{(-1)}_{+,+}$ (in view of Remark 1.1); here
$\Lambda ^{(-1)}_{+,+}\colon \dot H^0_p\simto \dot H^1_p$, but since
$\Lambda ^{(1)}_{-,+}\colon \ol H^1_p\simto \ol H^0_p$, it maps the
subspace $\dot H^1_p$ onto a subspace of $\ol H^0_p$ with infinite
codimension.
\endexample

Applications to fractional powers $A^a$ will be given below in Section
7.

\head 5. The $H_p^{\mu (s)}$-spaces and their boundary values \endhead

It will now be shown that the $H_p^{\mu (s)}$-spaces admit a special
definition of $\mu $-boundary values.

Let $M$ be a positive integer. First we consider $\E_\mu $ and
$\E_{\mu +M}$ for a smooth subset $\Omega $ of a paracompact manifold
$\Omega _1$ as in Section 2.

Let us
introduce the natural mapping
$$
\varrho _{\mu ,M} \colon \Cal E_\mu \to \Cal E_\mu /\Cal E_{\mu +M}.\tag5.1
$$
The first step is to represent $ \Cal E_\mu /\Cal E_{\mu +M} $ as the
space of sections of a trivial bundle and introduce norms in it. To do
so we first choose a Riemannian metric in $\Omega _1$ and then a $C^\infty $
function $d$ in $\comega$ which is equal to the distance from
$\partial\Omega $ sufficiently close to the boundary and is positive
and $C^\infty $ throughout $\Omega $. Set $$
I^\mu (x)=d(x)^\mu /\Gamma (\mu
+1)\text{ in }\comega\text{,  and }I^\mu =0\text{ in }\complement\Omega,\tag5.2$$ when
$\operatorname{Re}\mu >-1$ (consistently with (2.2)). 
This definition can be uniquely extended
modulo $C^\infty _0(\Omega )$ to arbitrary values of $\mu $ so that
$\partial_nI^\mu =I^{\mu -1}$, where $\partial_n$ denotes differentiation along the
geodesics perpendicular to $\partial\Omega $, sufficiently close to
$\partial\Omega $, and is defined as a $C^\infty $ function
elsewhere. By our definition of $\Cal E_\mu $ it follows easily that
every class in $ \Cal E_\mu /\Cal E_{\mu +1}$ contains an element of
the form $I^\mu (x)f$ where $f\in C^\infty (\comega )$, and
that such elements are congruent to 0 if and only if $f=0$ on the
boundary. By repeated application of this fact we conclude that any
element $u\in \Cal E_\mu $ can be written
$$
u=u_0I^\mu +u_1I^{\mu +1}+\dots+u_{M -1}I^{\mu +M-1}+v,\tag 5.3
$$
where the $u_j\in  C^\infty (\comega )$ are constant close to
$\partial\Omega $ on normal geodesics, and $v\in \Cal E_{\mu +M}$. The
boundary values of $u_j$ are uniquely determined by $u$, and it is
natural to write
$$
\gamma _{\mu ,j}u=u_j|_{\partial\Omega }.\tag5.4
$$
Note that
$$
\aligned
\gamma _{\mu ,j}u&=\gamma _{\mu +j,0}u, \text{ when }u\in \E_{\mu
+j};\\
\gamma _{\mu ,0}u&=\Gamma (\mu +1)\gamma _0 d(x)^{-\mu }u,\text{ when
$u\in\E_\mu $ with }\operatorname{Re}\mu >-1.
\endaligned\tag5.5
$$
When $\Omega =\rnp$, and $u(x)$ is written as $I^\mu w$ with
$I^\mu (x_n) =x_n^\mu /\Gamma (\mu +1)$ and $w(x)\in C^\infty (\crnp)$, then
$u_j(x')=\partial_n^jw(x',0)/\binom \mu  j$, where $\binom\mu  j=\Gamma
(\mu +j+1)/(j!\Gamma (\mu +1))$.

The mapping 
$$
\varrho _{\mu ,M}\colon u\mapsto \{\gamma _{\mu ,j}u\}_{j=0}^{M -1}\tag5.6
$$
has nullspace $\Cal E_{\mu +M}$ and identifies $ \Cal E_\mu /\Cal
E_{\mu +M}$ with $C^\infty (\partial\Omega )^{M }$; the mapping identifies with
the mapping in (5.1). The
identification depends of course on the choice of the Riemannian
structure but we shall keep it fixed in all that follows.
We can now think of $\varrho _{\mu ,M} $ as a mapping of
$\Cal E_\mu $ onto $C^\infty (\partial\Omega )^{M}$. 

\proclaim{Theorem 5.1} Let $s>\operatorname{Re}\mu +M-1/p'$, and let
$\comega $ equal $\crnp$ or a compact smooth manifold with boundary. The
mapping $\varrho _{\mu ,M}$ in {\rm (5.6)} extends by continuity to a
continuous mapping, also denoted $\varrho _{\mu ,M}$,
$$
\varrho _{\mu ,M}\colon H_p^{\mu (s)}(\comega)\to  \prod _{0\le j<M}B^{s-\operatorname{Re}\mu -j-1/p}_p(\partial\Omega );\tag5.7
$$
surjective and with kernel $H_p^{(\mu +M)(s)}(\comega)$. In other words, $\varrho
_{\mu ,M}$ defines a homeomorphism of $H_p^{\mu (s)}(\comega)/H_p^{(\mu
+M)(s)}(\comega)$ onto $ \prod _{0\le j<M}B^{s-\operatorname{Re}\mu -j-1/p}_p(\partial\Omega )$.
 \endproclaim

\demo{Proof}
We want to introduce in $ \Cal E_\mu /\Cal E_{\mu +M}$ the quotient
of the topology of $H_p^{\mu (s)}$. 
When discussing the quotient topology
it is 
sufficient to consider sections with support in a
local coordinate patch.

Thus let $u\in \Cal E_\mu (\crnp)\cap \Cal E'(K)$ where $K$ is a
compact set, and let $d(x)=x_n$. 
Writing $u$ in the form (5.3) we
have for $|\xi _n|>1$, say, and any $N$, 
$$
\hat u(\xi )=\sum_{j=0}^{M -1} b_j\hat u_j(\xi ')(\xi _n^-)^{-\mu
-j-1}+O([\xi ']
^{-N}|\xi _n|^{-\operatorname{Re}\mu -M-1}),\text{ where $b_j=i^{- (\mu +j+1)}$, }
$$
cf.\ (2.4).
This is similar to the formula (4.2), except that the nonzero factors
$b_j$ were incorporated in $\hat u_j$ in (4.2). 
Then we can use the calculation in (4.4) to obtain:
$$
\multline
\F(\Xi _+^\mu u)=i^\mu \hat u(\xi )(\xi _n-i[\xi '] )^\mu 
=i^\mu \sum_{j=0}^{M -1} \sum_{k=0}^jc_{jk}b_k\hat u_k(\xi ')[\xi ']
^{j-k}\xi _n^{ -j-1}+O([\xi ']^{-N}|\xi _n|^{-M-1})\\
= \sum_{j=0}^{M -1} \sum_{k=0}^jc_{jk}i^{ -k-1}\hat u_k(\xi ')[\xi ']
^{j-k}\xi _n^{ -j-1}+O([\xi ']^{-N}|\xi _n|^{-M-1}),
\endmultline
$$
where the $c_{jj}$ equal 1. Moreover, when $l<M$, 
$$
\F(\partial_n^l\Xi _+^\mu u)=(i\xi _n)^l\F(\Xi ^\mu _+u)
= \sum_{j=0}^{M -1} \sum_{k=0}^jc_{jk}i^{l -k-1}\hat u_k(\xi ')[\xi ']
^{j-k}\xi _n^{l -j-1}+O([\xi ']^{-N}|\xi _n|^{-2}).
$$
To calculate the boundary value $\gamma _0\partial_n^l\Xi ^\mu _+u$
from $\rnp$, note that for $l-j-1\ge 0$ the terms contribute with
distributions supported by $x_n=0$, and for $l-j-1<0$ it is the
coefficient of $\xi _n^{-1}$ that gives the boundary value at $x_n=0$,
cf.\ (2.9), so only $l=j$ contributes:
$$
\gamma _0\partial_n^j\Xi _+^\mu u
=
\gamma _0\Cal F^{-1}i\sum_{k=0}^j c_{jk}i^{j-k-1}\hat u_k(\xi ')[\xi ']
^{j-k}=\sum_{k=0}^j c'_{jk}[D']^{j-k} u_k, \tag 5.8
$$ 
with $c'_{jj}=1$ for all $j$.
In other words, with $\gamma _j=\gamma _0\partial_n^j$, the boundary values $\gamma _j \Xi ^\mu
_+u$ satisfy
$$
\pmatrix \gamma _0 \Xi _+^\mu u\\
\gamma _1 \Xi _+^\mu u\\
\vdots\\
\gamma _{M-1} \Xi _+^\mu u
\endpmatrix = \pmatrix 1& 0&\hdots&0\\ c'_{10}[D']& 1&\hdots &0\\
\vdots &\vdots& \ddots  &\vdots\\
c'_{M-1,0}[D']^{M-1}& c'_{M-1,1}[D']^{M-2}&\hdots&1
\endpmatrix \pmatrix \gamma _{\mu ,0} u\\\gamma _{\mu ,1}u\\ \vdots \\\gamma _{\mu ,M-1}u\endpmatrix
= \Phi \varrho _{\mu ,M} u,\tag5.9
$$
with an invertible triangular transition  matrix $\Phi $. 

Now we have from the well-known continuity 
properties of $\varrho _M=\{\gamma _0,\dots, \gamma _{M-1}\}$ (cf.\
(1.6)) that 
$$
\sum_{j=0}^{M -1}\|\gamma _j  \Xi ^\mu
_+u\|_{B^{s-\operatorname{Re}\mu -j-1/p}_p(\R^{n-1})}
\le C\|
r^+ \Xi ^\mu _+u\|^{\;}_{\ol H^{s-\operatorname{Re}\mu
}_p(\rnp)}
= C\|u\|_{\mu (s)}.
$$
Moreover, $\Phi $ is clearly a homeomorphism in $ \prod _{0\le
j<M}B^{s-\operatorname{Re}\mu -j-1/p}_p(\R^{n-1} )$, so by
(5.9), we likewise have
$$
\sum_{j=0}^{M -1}\|\gamma _{\mu ,j}u\|_{B^{s-\operatorname{Re}\mu -j-1/p}_p(\R^{n-1})}\le C
\|u\|_{\mu (s)}.\tag5.10
$$
Thus the mapping $\varrho _{\mu ,M}$ extends by continuity as
asserted. 

Finally, the extended map is surjective: For a given vector
$\varphi =\{\varphi _0,\dots,\varphi _{M-1}\}\in \prod _{0\le
j<M}B^{s-\operatorname{Re}\mu -j-1/p}_p(\R^{n-1})$, let $g\in
\ol H_p^{s-\operatorname{Re}\mu }(\rnp)$ be an element of $\ol
H^{s-\operatorname{Re}\mu }_p(\rnp)$ with $\varrho _Mg=\Phi \varphi $,
e.g.\ $g=\Cal K_M\Phi \varphi $ with $\Cal K_M$ defined in Section 1.1,
cf.\ (1.7). Set $u=\Xi ^{-\mu }_+e^+g$. By Proposition 1.7,
it has the desired properties. 
\qed
\enddemo

One can replace $[\xi ']$ by $\ang{\xi '}$ throughout the proof if convenient.

Note that on the space $H_p^{\mu (s)}(\comega)$, all the boundary
operators $\gamma _{\mu ,j}$, $j=0,1,\dots, M-1$, are defined when
$s>\operatorname{Re}\mu +M-1/p'$. They are {\it local}, in the sense
that they are extensions by continuity of local operators of the form: $\gamma _0$ composed with
multiplication and differential operators.
For this extended definition, the first line in (5.5) is valid on $H^{(\mu +j)(s)}(\comega)$,
and the second line holds on $H^{\mu (s)}(\comega)$ when
$\operatorname{Re}\mu >-1$.

\example{Remark 5.2} In the course of the above proof we have in fact
constructed an explicit right inverse to
$\varrho _{\mu ,M}$ in the case $\Omega =\rnp$, namely 
$$
\Cal K_{\mu ,M}=\Xi ^{-\mu }_+e^+\Cal K_M\Phi .\tag 5.11
$$
\endexample

We observe in particular from (5.9) that $\Phi =I$ when $M=1$, and hence $\gamma _0  \Xi _+^\mu u
=\gamma _{\mu ,0} u$. For the case $M=1$ we consequently have:

\proclaim{Corollary 5.3} When $s>\operatorname{Re}\mu +1/p$, the
mapping $\gamma _{\mu ,0} $ 
is continuous and surjective from $H_p^{\mu (s)}(\crnp)$ 
to
$B^{s-\operatorname{Re}\mu -1/p}_p(\R^{n-1} )$ 
with nullspace $H_p^{(\mu
+1)(s)}(\crnp)$.
It coincides
with $ \gamma _0\Xi _+^\mu $. A right inverse is $K_{\mu
,0}=\Xi ^{-\mu }_+e^+K_0$, where $K_0\colon B^{t
-1/p}_p(\R^{n-1} )\to \ol H_p^{t}(\rnp)$ is a right inverse
of $\gamma _0$.
\endproclaim

As an example, let us also do the calculation of $\Phi $ in detail in the case $M=2$.

For
$u\in \E_\mu (\crnp)\cap \E'(K)$,
$$
u(x',x_n)=u_0(x')I^\mu (x_n)+u_1(x')I^{\mu
+1}(x_n)+\text{remainder},
$$
so we have for $|\xi _n|\ge 1$ (assumed in the following):
$$
\hat u(\xi )=i^{-\mu -1}\hat u_0(\xi ')(\xi _n^-)^{-\mu -1}+i^{-\mu
-2}\hat u_1(\xi ')(\xi _n^-)^{-\mu -2} +O(\xi _n^{-\mu -3}).
$$
Denote $[\xi ']=\sigma $. The function $(\sigma +i\xi _n)^\mu $ is
Taylor expanded:
$$
(\sigma +i\xi _n)^{\sigma }=i^\mu (\xi _n-i\sigma )^\mu =i^\mu (\xi
_n^-)^\mu 
-i^{\mu -1}\mu \sigma (\xi _n^-)^{\mu -1}+O(\xi _n^{\mu -2}). 
$$
Hence 
$$
(\sigma +i\xi _n)^{\sigma }\hat u(\xi )=i^{-1}\hat u_0(\xi ') \xi _n^{-1}
+i^{ -2}\mu \sigma \hat u_0(\xi ')\xi _n^{ -2}
+i^{ -2} \hat u_1(\xi ')\xi _n^{
-2}
+O(\xi _n^{-3}). 
$$
In view of (2.9),
$$
\gamma _0\Xi _+^\mu u=u_0.
$$
Moreover,
$$
i\xi _n(\sigma +i\xi _n)^{\sigma }\hat u(\xi )=\hat u_0(\xi ') 
+i^{ -1}\mu \sigma \hat u_0(\xi ')\xi _n^{ -1}
+i^{ -1} \hat u_1(\xi ')\xi _n^{
-1}
+O(\xi _n^{-2}), 
$$
 so since $\F^{-1}_{\xi \to x}\hat u_0(\xi ')=u_0(x')\otimes \delta
 _0(x_n)$ does not contribute to the boundary value from $\rnp$,
$$
\gamma _0\partial_n\Xi ^\mu _+u =\mu \,\sigma (D')\,u_0 +u_1.
$$
Thus
$$
\pmatrix \gamma _0\Xi ^\mu _+u\\ \gamma _1\Xi ^\mu
_+u\endpmatrix=\pmatrix 1&0\\ \mu \,[D']&1
\endpmatrix \pmatrix \gamma _{\mu ,0}u\\ \gamma _{\mu ,1}u\endpmatrix,
\text{ and }\Phi =\pmatrix 1&0\\ \mu \,[D']&1\endpmatrix.\tag5.12
$$
If $\sigma $ is taken  equal to $\ang{\xi '}$ instead of  $[\xi ']$, we get of
course $\Phi $ of the above form with $[D']$ replaced by $\ang{D'}$.
\medskip

By use of concrete formulas from the Boutet de Monvel calculus we can
show that not only the boundary operators from $H_p^{\mu (s)}(\crnp)$
carry a $\mu $'th power of $x_n$, but also the functions on $\rnp$
themselves do so.

\proclaim{Theorem 5.4} When $s>\operatorname{Re}\mu +M-1/p'$, and
$u\in H_p^{\mu (s)}(\crnp)$, then with $\Cal K_{\mu ,M}$ taken as in
{\rm (5.11)},
$$
u=v+w,\text{ where }v=\Cal K_{\mu ,M}\varrho _{\mu ,M}u\text{ and } 
w\in H^{(\mu +M)(s)}(\crnp).\tag5.13
$$
Here if $\operatorname{Re}\mu >-1$,  $v=\Xi _+^{-\mu }e^+\Cal K_M \varrho _M\Xi _+^\mu u$ has the form 
$$
v
={\sum}_{j=0}^{M-1}c_jx_n^{\mu +j}e^+K_0(\gamma _{\mu ,j}u)=e^+x_n^\mu
v_0,
\tag5.14 
$$
with $ v_0\in \ol H^{s-\operatorname{Re}\mu }(\rnp)$, $K_0$ as in {\rm
(1.7)}.

Thus one has for
$\operatorname{Re}\mu >-1$, $s>\operatorname{Re}\mu -1/p'$, with $M\in{\Bbb N}$:
$$
\aligned
&H^{\mu (s)}_p(\crnp)\cases =\dot H^s_p(\crnp)\text{ if
}s-\operatorname{Re}\mu \in \,]-1/p',1/p[\,,\\
\subset\dot H^{s-0}_p(\crnp)\text{ if }s-\operatorname{Re}\mu =1/p.\endcases
\\
H^{\mu (s)}_p(\crnp)&\subset e^+x_n^{\mu }\ol H_p^{s-\operatorname{Re}\mu
}(\rnp)+
\cases \dot H^s_p(\crnp)\text{ if }s-\operatorname{Re}\mu \in M+\,]-1/p',1/p[\, \\
\dot H^{s-0}_p(\crnp)\text{ if }s-\operatorname{Re}\mu = M+1/p .\endcases
\endaligned\tag 5.15
 $$

The inclusions {\rm (5.15)} also hold in the manifold situation, with
$\rnp$ replaced by $\Omega $ and $x_n$ replaced by $d(x)$.

\endproclaim

\demo{Proof} 
The decomposition (5.13) is an immediate consequence of Theorem 5.1;
here $w\in H_p^{(\mu +M)(s)}(\crnp)$ since $\varrho _{\mu ,M}w=0$. 
In the next statements we take $\operatorname{Re}\mu >-1$ in order to identify $I^\mu $ with
the locally integrable function $e^+r^+x_n^\mu /\Gamma (\mu +1)$. Distributional
formulations can be made for lower $\mu $.

For the  description in (5.14), note that the first equality follows
from (5.9) and (5.11). For the next equality,  consider first the case $M=1$,
where simply $v=
K_{\mu ,0}\gamma _{\mu ,0}u$.

Recall from (1.7) that $K_0$ is the elementary Poisson operator of
order 0 $$
\varphi \mapsto \F^{-1}_{\xi
'\to x'}\bigl(\hat\varphi (\xi ')e^+r^+e^{-[\xi ']x_n}\bigr)=\F^{-1}_{\xi
\to x}\bigl(\hat\varphi (\xi ')([\xi ']+i\xi _n)^{-1}\bigr).
$$ Constructing $K_{\mu,0}$ as in Corollary 5.3 we have, cf.\  (2.5),
$$
\aligned
K_{\mu,0}\varphi &= \F^{-1}_{\xi \to x}\bigl(([\xi
']+i\xi _n)^{-\mu }\hat\varphi (\xi ') ([\xi ']+i\xi _n)^{-1}\bigr)\\
&=c_\mu\F^{-1}_{\xi '\to x'}\bigl(e^+r^+x_n^{\mu}e^{-[\xi ']x_n}
\hat\varphi (\xi ')\bigr)
=c_\mu e^+x_n^{\mu}K_0\varphi .
\endaligned\tag 5.16
$$
Hence since $\gamma _{\mu,0}u\in B_p^{s-\operatorname{Re}\mu -1/p}({\Bbb R}^{n-1})$,
$$
v=c_\mu e^+x_n^\mu K_0\gamma _{\mu,0}u\in e^+x_n^\mu\ol H_p^{s-\operatorname{Re}\mu }(\rnp)
,\tag5.17
$$
 by the mapping properties of Poisson operators shown in \cite{G90}.

For general $M$ we have that $v=K_{\mu ,0}\gamma _{\mu ,0}u+\dots
+K_{\mu ,M-1}\gamma _{\mu ,M-1}u$, and we have to account for the
general term $K_{\mu ,j}\gamma _{\mu ,j}u$. 
Here $\varphi _j=\gamma
 _{\mu ,j}u\in B_p^{s-\operatorname{Re}\mu -j-1/p}({\Bbb R}^{n-1})$.
By (1.7), $K_j$ acts as 
$$
\varphi _j \mapsto \F^{-1}_{\xi
\to x}\bigl(\hat\varphi _j(\xi ')\tfrac{(-1)^j}{j!}\partial_{\xi _n}^j([\xi ']+i\xi _n)^{-1}\bigr)=\F^{-1}_{\xi
\to x}\bigl(\hat\varphi _j(\xi ')i^j([\xi ']+i\xi _n)^{-j-1}\bigr).
$$ Then
$$
\aligned
K_{\mu,j}\varphi _j&= \F^{-1}_{\xi \to x}\bigl(([\xi
']+i\xi _n)^{-\mu }\hat\varphi _j(\xi ') ([\xi ']+i\xi _n)^{-j-1}\bigr)\\
&=c_{\mu ,j}\F^{-1}_{\xi '\to x'}\bigl(e^+r^+x_n^{\mu+j}e^{-[\xi ']x_n}
\hat\varphi _j(\xi ')\bigr)
=c_{\mu ,j}e^+x_n^{\mu+j}K_0\varphi _j.
\endaligned\tag5.18
$$
By the rules of the Boutet de Monvel calculus, $x_n^jK_0$ is
a Poisson operator of order $-j$, so the mapping properties from
\cite{G90} assure that $x_n^jK_0\varphi _j\in \ol
H_p^{s-\operatorname{Re}\mu }(\rnp)$. Thus
$$
K_{\mu ,j}\gamma _{\mu ,j}u\in e^+x_n^\mu \ol
H_p^{s-\operatorname{Re}\mu }(\rnp).
$$

The first line in (5.15) is shown in (1.26) when
$s-\operatorname{Re}\mu <1/p$, and when $s-\operatorname{Re}\mu =1/p$,
it follows in view of (1.31).
The second line in (5.15) follows from (5.13) and (5.14), when $s-
\operatorname{Re}\mu -M\in \,]-1/p',1/p[\,$, since $H^{(\mu
+M)(s)}(\crnp)$ then is as in the first line.

The conclusions in (5.15) carry over to the manifold situation by
use of local coordinates. 
\qed

\enddemo

The formulas (5.17), (5.18) are of interest in themselves.

\proclaim{Corollary 5.5} Let $\operatorname{Re}\mu \ge 0$,
$s>\operatorname{Re}\mu +n/p$. Then 
$$
H_p^{\mu (s)}(\comega)\subset e^+d(x)^\mu  C^{s-\operatorname{Re}\mu -n/p-0}(\comega),\tag5.19
$$
where $-0$ can be left out when $s-\operatorname{Re}\mu -n/p$,
$s-n/p$ and $s-\operatorname{Re}\mu -1/p$ are noninteger.
\endproclaim

\demo{Proof} We use the description by two terms in (5.15). By (1.23),
$$
e^+d(x)^\mu \ol H_p^{s-\operatorname{Re}\mu }(\Omega)\subset e^+d(x)^\mu  C^{s-\operatorname{Re}\mu -n/p-0}(\comega),
$$
where $-0$ can be left out when $s-\operatorname{Re}\mu -n/p$ is not
integer.
When $u\in \dot H_p^s(\comega)$, it belongs to $C^{s-n/p-0}(\Omega
_1)$ and is supported in $\comega$; here $-0$ can be left out when $s-n/p$
is not integer. Since $s>1/p$, $\gamma _0u=0$;
then 
in view of the H\"older continuity, $u\in e^+d(x)^\mu C^{s-\operatorname{Re}\mu -n/p-0}(\comega)$,
since $s-n/p>\operatorname{Re}\mu \ge 0$. This extends to $\dot
H_p^{s-0}(\comega)$ when $s-\operatorname{Re}\mu -1/p$ is integer; the $-0$ is needed then
in view of (5.15). Hereby the assertion is verified for the two terms in
(5.15).\qed\enddemo

\head 6. Nonhomogeneous boundary value problems, parametrices \endhead

The problems treated in Theorem 4.4 can be regarded as homogeneous boundary
problems, when we see them in the following perspective.

Consider again our operator $P$ satisfying the hypotheses of Theorem 4.4, with the
factorization index
 $\mu _0\in{\Bbb C}$. For a positive integer $M$ let $\mu =\mu _0-M$.
We have from Theorem 5.1 that when $s>\operatorname{Re}\mu
+M-1/p=\operatorname{Re}\mu _0-1/p$, then $\varrho _{\mu ,M}$ defines a homeomorphism
$$
\varrho _{\mu ,M}\colon H_p^{\mu (s)}(\comega)/H_p^{\mu _0(s)}(\comega)
\simto  \prod _{0\le j<M}B^{s-\operatorname{Re}\mu -j-1/p}_p(\partial\Omega ).\tag6.1
$$
Combining this  with the Fredholm property of 
$$
r^+P\colon H_p^{\mu _0(s)}(\comega)\to \ol H_p^{s-\operatorname{Re}m}(\Omega ),\tag6.2
$$
we have immediately:

\proclaim{Theorem 6.1} Let $P$ satisfy the hypotheses of Theorem {\rm 4.4},
and let $\mu =\mu _0-M$ for a positive integer $M$. Then when $s>\operatorname{Re}\mu
_0-1/p'$,  $\{r^+P,\varrho _{\mu ,M}\}$ 
defines a Fredholm operator
$$
\{r^+P,\varrho _{\mu ,M}\}\colon H_p^{\mu (s)}(\comega)\to \ol
H_p^{s-\operatorname{Re}m}(\Omega )\times \prod _{0\le j<M}B^{s-\operatorname{Re}\mu -j-1/p}_p(\partial\Omega ).\tag6.3
$$

\endproclaim

This is a solvability result for the following inhomogeneous ``Dirichlet problem'' for $P$:
$$
r^+Pu=f,\quad \varrho _{\mu ,M}u=\varphi ,\tag6.4
$$
where $\varphi $ is an $M$-vector $\{\varphi _0,\dots, \varphi
_{M-1}\}$ of boundary data. 

We can in particular take $M=1$; this gives:

\proclaim{Corollary 6.2} 
With $P$ as in Theorem {\rm 5.1}, let
$\mu =\mu _0-1$. Then  
$$
\{r^+P,\gamma _{\mu ,0}\}\colon H_p^{\mu (s)}(\comega)\to \ol
H_p^{s-\operatorname{Re}m}(\Omega )\times B^{s-\operatorname{Re}\mu -1/p}_p(\partial\Omega )\tag6.5
$$
is Fredholm when $s>\operatorname{Re}\mu +1-1/p'(=
\operatorname{Re}\mu _0-1/p')$. 
\endproclaim

This shows a solvability result for the problem 
$$
r^+Pu=f,\quad \gamma _{\mu ,0}u=\varphi _0.\tag6.6
$$
with just $\gamma _{\mu ,0}u$
prescribed, $\mu =\mu _0-1$.

\example{Example 6.3}
For the Laplace-Beltrami operator, $\mu _0=1$, so Corollary 6.2 is applicable with $\mu =0$. Here
$H^{0(s)}_p=\ol H^s_p$ and $\gamma _{0,0}= \gamma _0$, so it gives
the Fredholm property of the mapping
$$
\{\Delta ,\gamma _0\}\colon \ol H^s_p(\Omega )\to \ol H^{s-2}_p(\Omega )\times B^{s-1/p}_p(\partial\Omega )
$$
for $s>1/p$, which is well-known as the inhomogeneous Dirchlet problem
for $\Delta $.

For $M=2$, $\mu =\mu _0-M=-1$ and $\varrho _{\mu ,M}=\{\gamma _{-1,0},
\gamma _{-1,1}\}$. When $u\in \E_{-1}(\crnp)$,
$$
u=u_0(x')\delta (x_n)+u_1(x')+v,\quad v\in \E_1(\crnp), u_0\text{ and }u_1\in C^\infty (\R^{n-1}),
$$
according to (5.3); then $\gamma _{-1,0}u=u_0(x')$ and $\gamma
_{-1,1}u=u_1(x')$. We get a solvability result for $\Delta $ where the term 
$u_0(x')\delta (x_n)$ can be prescribed arbitrarily. This is a point
of view on boundary problems related to the works of Roitberg and
Sheftel' \cite{RS69}, \cite{R96}, going beyond the ordinary concept of boundary value problems.
\endexample

\example{Remark 6.4} Since the distributions $I^\mu (x_n)$ are
locally integrable functions $e^+r^+c_\mu x_n^\mu $ only when
$\operatorname{Re}\mu >-1$, 
the trace maps $\gamma _{\mu ,0}$ are
somewhat ``wild'' when $\operatorname{Re}\mu \le -1$. In the
interpretations of concrete cases we shall in this paper only consider
situations where the entering trace operators have
$\operatorname{Re}\mu >-1$; e.g.\ in applications of Theorem 6.1 we
only take $M<\operatorname{Re}\mu _0+1$.
\endexample

We shall finally show that a parametrix of the 
nonhomogeneous boun\-da\-ry problem considered in 
Corollary 6.2 can be obtained by a combination of the knowledge from
the type 0 calculus and the special operators used here. The
construction of $K$ ``from scratch'' takes up much effort in \cite{H65}.

\proclaim{Theorem 6.5} Let $P$ be a globally estimated $\psi $do of
order $m\in{\Bbb C}$ and type $\mu _0\in{\Bbb C}$, and factorization
index $\mu _0$, relative to the domain $\crnp$. Let $s>\operatorname{Re}\mu _0-1/p'$.

For the problem considered in Corollary {\rm 6.2}:
$$
r^+Pu=f,\quad \gamma _{\mu _0-1,0}u=\varphi ,
\tag6.7$$
with $f$ given in $\ol H_p^{s-\operatorname{Re}m}(\rnp )$ and $\varphi $ given in $B^{s-\mu
_0+1-1/p}_p(\R^{n-1} )$, a parametrix is 
$$
\pmatrix R & K\endpmatrix \colon \matrix \ol
H^{s-\operatorname{Re}m}_p(\rnp )\\ 
\times \\ B^{s-\mu
_0+1-1/p}_p(\R^{n-1} ) \endmatrix \to H^{(\mu _0-1)(s)}_p(\crnp ) ,  \tag6.8
$$
 where $R$ is as in Theorem {\rm 4.4}, and $K$ is of the form
$$
K 
=\Xi _+^{1-\mu _0}e^+K' =\Lambda _+^{1-\mu _0}e^+K'' ,
\tag 6.9
$$
 with  Poisson operators $K'$ and $K''$ of order $0$ in
the Boutet de Monvel calculus.

\endproclaim 

\demo{Proof}
As a parametrix for the problem (6.7) with $\varphi =0$ we can use  $R$
introduced in Theorem 4.4,
since $H^{\mu _0(s)}_p$ is the subspace of $H^{(\mu _0-1)(s)}_p$ where
$\gamma _{\mu _0-1,0}u=0$. Note that $P$ is expressed in terms of $Q$ by
$$
P
=\Lambda ^{m-\mu _0}_{-}   Q \Lambda ^{\mu _0}_+ .\tag6.10
$$

It remains to solve problem (6.7) 
when $f=0$. Consider  
$$
r^+Pu=0,\quad \gamma _{\mu _0-1,0}u=\varphi ,
\tag6.11$$
with $\varphi $ given in $B^{s-\mu _0+1-1/p}_p(\R^{n-1} )$. 
On $\rnp$ we have explicit formulas for
the elementary Poisson-like operators $\Cal K_{\mu ,M}$. Here 
$$
K_{\mu _0-1,0}=\Xi _{+}^{1-\mu _0}e^+K_0,\tag6.12
$$
  cf.\ Corollary 5.3. To solve (6.11), let 
$$
z = \Xi _{+}^{1-\mu _0}e^+K_0\varphi ,
$$
and form $w =u-z $; it must solve
$$
r^+Pw =-r^+P\Xi _{+}^{1-\mu _0}e^+K_0\varphi ,\quad \gamma _{\mu _0-1,0}w =0.
\tag6.13
$$
By Theorem 4.4, this problem has the solution in a parametrix sense:
$$
\aligned
w &=-Rr^+P\Xi _{+}^{1-\mu _0}e^+K_0\varphi 
=-\Lambda ^{-\mu
 _0}_+e^+\widetilde{Q_+}\Lambda ^{\mu _0-m}_{-,+}r^+\Lambda ^{m-\mu _0}_{-}   Q \Lambda ^{\mu _0}_+ \Xi _{+}^{1-\mu _0}e^+K_0\varphi 
\\
&=-\Lambda ^{-\mu
 _0}_+e^+\widetilde{Q_+} r^+ Q \Lambda ^{\mu _0}_+ \Xi _{+}^{1-\mu _0}e^+K_0\varphi , \endaligned
$$
when we take (6.10) into account, using also Remark 1.1.

We now observe, recalling the definition of $Y^\mu _+$ from (1.16)ff., that
$$
 \Lambda ^{\mu _0}_+ \Xi _{+}^{1-\mu _0}e^+K_0= \Lambda
 _+^1e^+r^+\operatorname{OP}(\lambda _+^{\mu _0-1}\chi _+^{1-\mu _0})e^+K_0=
\Lambda
 _+^1e^+Y^{\mu _0-1}_{+,+}K_0.
$$
An application of Lemma 6.6 below gives that $Y^{\mu _0-1}_{+,+}K_0$
is a Poisson operator of order 0 in the Boutet de Monvel calculus.
Hence
$\widetilde{Q_+} r^+ Q \Lambda ^{1}_+ e^+Y_{+,+}^{\mu _0-1}K_0$
is a Poisson operator $K_1$ of order 1 in the Boutet de Monvel
calculus, and
$$
w=-\Lambda ^{-\mu _0}_+e^+K_1\varphi .\tag6.14
$$
This can be rewritten, using again Lemma 6.6, as 
$$
w=
-\Xi ^{-\mu _0}_+e^+Y^{-\mu _0}_{+,+}K_1\varphi ,
$$
where $Y^{-\mu _0}_{+,+}K_1$ is another Poisson operator of order 1.
Thus $u=z+w$ has the structure
$$
u=\Xi _{+}^{1-\mu _0}e^+K'\varphi,
$$
with a Poisson operator $K'$ of order 1. This shows the first formula
in (6.9). For the second formula, we keep $w$ in the form (6.14) and
instead rewrite 
$$
z=\Xi _{+}^{1-\mu _0}e^+K_0\varphi=\Lambda  _{+}^{-\mu _0}e^+\Lambda ^1_{+,+}Y^{\mu
_0-1}_{+,+}K_0\varphi ,
$$
where $\Lambda ^1_{+,+}Y^{\mu
_0-1}_{+,+}K_0$ is a Poisson operator of order 1 in view of Lemma 6.6
and the composition rules.
\qed

\enddemo

 Analogous constructions can be made in case $M>1$.

The following lemma shows a case where that the composition of a
Poisson operator with certain generalized $\psi $do's defined from
symbols that only satisfy some of the estimates required for the
$S^d_{1,0}$ classes, is again a Poisson operator (similarly to some cases considered in Sect.\
3.2 of \cite{GK93}).

\proclaim{Lemma 6.6} Let $K$ be a Poisson operator on $\crnp$ of order
$m$, with symbol $k(x',\xi )$, let $s(x',\xi )$ be a Poisson symbol of
order $0$, and let $S=\operatorname{OP}(s(x',\xi ))$ be the
generalized $\psi $do with symbol $s$ (defined as in {\rm (1.2)}). Then
the composed operator $S_+K$ is a Poisson operator of order $0$ with
symbol $k'= s \circ k\sim \sum_{\alpha \in {\Bbb N}_0^n}\frac
{1}{\alpha !}D_\xi ^\alpha s\, \partial_x^\alpha k$.

The result applies in particular when $S=Y^{\mu
}_+-1=\operatorname{OP}(\eta ^\mu _+)-1$ as defined in {\rm (1.6)ff.}
\endproclaim

\demo{Proof} When $k$ is independent of $x'$, so that $e^+Ku=\Cal
F^{-1}(k(\xi )\hat u(\xi '))$, we can move $k(\xi )\hat u$ inside the
integral defining the action of $S$, and the result follows since $sk$
is a Poisson symbol of order $m$ (a product of functions in $\Cal H^+$
is in $\Cal H^+$). In the
$x'$-dependent case, there is a standard procedure of replacing $k$ by
a $y'$-form symbol; it can then be
moved inside the integral as above, and the resulting symbol
in  $(x',y')$-form reduced to $x'$-form as an asymptotic series.

For the last statement, we recall that 
$$
\lambda ^1_+/\chi ^1_+=1+q^1_+(\xi ),\quad q^1_+=[\xi '](\bar\psi (\xi
_n/a[\xi '])-1)/([\xi ']+i\xi _n)\in \Cal H^+
$$
as a function of $\xi _n$ for all $\xi '$, and $|q^1_+(\xi )|\le
\frac12$ (recall that $a$ is taken large). Then
$$
\eta ^\mu _+(\xi )=(1+q^1_+)^\mu =1+\mu q^1_++\mu (\mu -1)\tfrac12 (q^1_+)^2+\dots=1+q
$$
as a convergent Taylor series, where $(q^1_+)^k\in \Cal H^+_{-k}$ as a
function of $\xi _n$, so that
$q\in \Cal H^+$ for all $\xi '$; moreover, it is homogeneous of degree 0 for $|\xi '|\ge 1$.  
\qed\enddemo

\proclaim{Corollary 6.7} The operator $K$ in Theorem {\rm 6.5} has the
property that
when $B$ is a $\psi $do of type $\mu _0$ and order $m_0+\mu _0$,  $m_0\in{\Bbb Z}$,
then $\gamma _0r^+BK$ is a $\psi $do on ${\Bbb R}^{n-1}$ of order $m_0+1$.
\endproclaim

\demo{Proof} Let $B'=B\Lambda _+^{1-\mu _0}$; then $B'$ is of order
$m_0+1$ and type 0, hence belongs to the Boutet de Monvel
calculus. From the rules there we conclude, using (6.9), 
that $\gamma _0r^+BK=\gamma _0B'_+K''$
is a $\psi $do of order $m_0+1$.\qed
\enddemo

\head 7. Applications to fractional powers of elliptic operators \endhead

We here show some consequences for fractional
powers of differential operators.
Let $A$ be a second-order strongly elliptic operator with $C^\infty
$-coefficients on $\Omega _1$ (that can be taken compact), and consider the
fractional powers $P_a=A^{a}$ for $a>0$. By Lemma 2.9 and Example 3.2, they
are classical $\psi $do's of order $2a$, having type $a$ and 
factorization index $\mu _0=a$ relative to  $\Omega $. This holds in
particular for $(-\Delta )^a$, where $\Delta $ is the Laplace-Beltrami
operator on $\Omega _1$. See also Remark 2.10. 

We have as an immediate corollary of Theorems 4.4 and 6.1:

\proclaim{Theorem 7.1} Let $1<p<\infty $, and let $s>a-1/p'=a-1+1/p$.

$1^\circ$ Let $u\in \dot H_p^\sigma (\comega) $ for
some $\sigma >a-1/p'$. If $r^+ P_au\in \ol H_p^{s-2a}(\Omega )$,
then $u\in H_p^{a(s)}(\comega)$. The mapping $r^+P_a$  is Fredholm:
$$
r^+  P_a\colon H^{a(s)}_p(\comega )\to \overline
H^{s-2a}_p(\Omega ).\tag 7.1
$$

$2^\circ$ In particular, if $r^+ P_au\in C^\infty (\comega)$, then
$u\in \E_{a}(\comega)$, and the mapping  $r^+P_a$ is Fredholm: 
$$
r^+  P_a\colon\E_{a}(\comega)\to  C^\infty (\comega).\tag7.2
$$

$3^\circ$ Moreover, when $M$ is a positive integer, the operator
$\{r^+P_a,\varrho _{a-M ,M}\}$ is Fredholm:
$$
\aligned
\{r^+P_a,\varrho _{a-M ,M}\}
 &\colon H_p^{(a-M) (s)}(\comega)\to \ol
H_p^{s-2a}(\Omega )\times \prod _{0\le j<M}B^{s-a
+M-j-1/p}_p(\partial\Omega ),\\
\{r^+P_a,\varrho _{a-M ,M}\}
&\colon \E_{a-M}(\comega)\to C^\infty (\comega )\times C^\infty (\partial\Omega )^M.\endaligned\tag7.3
$$
\endproclaim

As mentioned in Remark 6.4, we shall here only discuss $3^\circ $
when $M<a+1$.

\example{Example 7.2} Let us describe the domain of the Dirichlet
realization for $p=2$ in this context. Define it as the space of solutions of
$r^+P_af=u$ with $f\in L_2(\Omega )$ according to the above theorem:
$$
D(P_{a,\operatorname{Dir}})=\{u\in \dot H_2^{a-\frac12+0}(\comega)\mid r^+P_au\in L_2(\Omega )\}.
$$
The order of $P_a$ is $2a$, so the range space in Theorem 7.1 $1^\circ$ equals $L_2(\Omega )$ when
$s=2a$. Then
$D(P_{a,\operatorname{Dir}})=H_2^{a(2a)}(\comega)$, where $r^+P_a$ is
Fredholm. This is a precise
and seemingly new result when $a\ge \frac12$, the case $a<\frac12$
being covered by Vishik and Eskin's theorem.

Note that$$
2a\in a+\,]-\tfrac12 ,\tfrac12 [\,\text{ when }a<\tfrac12 ,\; 2a\in
a+1+[-\tfrac12 ,\tfrac12 [\,\text{ when }\tfrac12 \le a<\tfrac 32,\text{ etc.}
$$
Then we have by Theorem 5.4,
$$
D(P_{a,\operatorname{Dir}})\cases = \dot H_2^{2a}(\comega), \text{ when }0<a<\tfrac12 ,\\
=H_2^{\frac12 (1)}(\comega)\subset \dot H_2^{1-0}(\comega)\text{ when
}a=\tfrac12 ,\\
\subset e^+d(x)^a\ol H_2^a(\Omega )+\dot H_2^{2a}(\comega)\text{ when
}\tfrac12 <a<\tfrac32,\text{ etc.}
\endcases\tag7.4
$$
For $a> \frac12$, the structure of the contribution from $d(x)^a\ol H^a_2$
is described in (5.14), (5.17).
\endexample

We remark that the operator $P_{a,\operatorname{Dir}}$ for $A=-\Delta $ is not
the same as the operator $B_a=(-\Delta _{\operatorname{Dir}})^a$
defined by $L_2$ spectral theory from the Dirichlet realization $\Delta _{\operatorname{Dir}}$ of the
Laplacian when $0<a<1$. Here $D(B_a)$ is the interpolation space
between $\ol H^2_2(\Omega )\cap \dot H^1_2(\comega)$ and $L_2(\Omega )$, equal to $\{u\in \ol
H_2^{2a}(\Omega )\mid \gamma _0u=0\}$ when $a>\frac14$ and to $\dot
H^{2a}_2(\comega)$ when $a<\frac34$.

Now we want to see what the result gives in terms of bounded or H\"o{}lder continuous
functions.  It has been shown by Ros-Oton and Serra in \cite{RS14} for
$0<a<1$, $\Omega \subset{\Bbb R}^n$,  that solutions of 
$r^+(-\Delta )^au=f\in L_{\infty }(\Omega )$ with
$u\in \dot H^a(\comega)$ 
are
in $d(x)^aC^\alpha (\comega)$ for some $\alpha <\min\{a,1-a\}$, when
$\Omega $ is $C^{1,1}$. (See \cite{RS14} for further references to
contributions to the problem.)

Let us study the solutions of the homogeneous Dirichlet problem 
$$
r^+ P_au=f,\tag7.5
$$ 
where $f$ is given in $\ol
H^t_p(\Omega )$ with $t\ge 0$,  for $u\in \dot H_p^{a-1/p'+0}(\comega)$.
By Theorem 7.1 $1^\circ$ with $s=t+2a$,  $u$ belongs to $H_p^{a(t+2a)}(\comega
)$. 
By Corollary 5.5,
$$
H_p^{a(t+2a)}(\comega
)\subset e^+d(x)^aC^{t+a-n/p-0}(\comega), \tag7.6
$$
when $p$ is so large that $a>n/p$ (for then $t+2a>a+n/p$); here
$-0$ can be left out except at certain values of $t$.
The ellipticity of $P_a$ moreover assures that $u\in H^{t+2a}_{p,\operatorname{loc}}(\Omega )$, which is contained
in $C^{t+2a-n/p -0}(\Omega )$. We conclude that
$$
u\in e^+ d(x)^a C^{t+a-n/p -0}(\comega )\cap C^{t+2a-n/p -0}(\Omega ).\tag7.7
$$

Note that the prerequisite $u\in \dot H_p^{a-1/p'+0}(\comega)$ is
satisfied if (cf.\ (1.23))
$$
u\in \cases e^+L_p(\Omega ),\text{ when }a<1/p',\\
\dot C^{a-1/p'+0}(\comega ),
\text{ when }a\ge 1/p'.\endcases \tag7.8
$$

For $t=0$ we have found in particular: 
$$
f\in L_p (\Omega ) \implies u\in e^+d(x)^a
C^{a-n/p-0 }(\comega)\cap C^{2a-n/p-0 }(\Omega ),\tag7.9
$$
where $-0$ can be omitted when $a-n/p$, $2a-n/p$ and $a-1/p$ are not integer.
For $p\to\infty $, $a-n/p\to a$, and (7.9) gives, since $L_\infty
(\Omega )\subset L_p(\Omega )$ for all $p$,
$$
f\in L_\infty (\Omega ) \implies u\in e^+d(x)^a 
C^{a-0 }(\comega)\cap C^{2a-0 }(\Omega ).\tag7.10
$$
(It suffices that  $u\in \dot H_{p_0}^{a-1/p_0'+0}(\comega)$ for some $p_0$.)

This shows an improvement of Th.\ 1.2 of Ros-Oton and Serra \cite{RS14}, in higher
generality concerning the studied operator and the data, when the
boundary is smooth. 

For general higher $t$, we similarly find, noting that
$C^{t+0}(\comega)\subset \ol H^{t}_p(\Omega )$ and letting $p\to\infty $:
$$
f\in  C^{t+0} (\comega )\implies u\in e^+d(x)^a C^{t+a-0
}(\comega )\cap C^{t+2a-0}(\Omega ).\tag7.11
$$

Recall also that Theorem 7.1 $2^\circ$ shows:
$$
f\in  C^{\infty } (\comega )\iff u\in e^+d(x)^a C^{\infty 
}(\comega )\,\bigl(=\E_{a}(\comega)\bigr),\tag7.12
$$
with Fredholm solvability, when  $u\in \dot
H_p^{a-1/p'+0}(\comega)$for some $p$.
 
This
extends results of \cite{RS14} to arbitrarily smooth spaces.
The Fredholm property of (7.1) implies  that in each of the cases (7.9)--(7.11), there is solvability
for $f$ in the indicated space, subject to a finite dimensional linear condition.

We have hereby obtained:

\proclaim{Theorem 7.3} Let $A$ be a second-order strongly elliptic
differential operator on $\Omega _1$ with smooth coefficients, and let
$P_a=A^a$ for some $a>0$, a $\psi $do of order $2a$ by Seeley's
construction. Let $d(x)>0$ on $\Omega $, $d\in C^\infty (\comega)$ and
proportional to 
$\operatorname{dist}(x,\partial\Omega )$ near $\partial\Omega $. Consider the
homogeneous Dirichlet problem {\rm (7.5)}, taking $u\in \dot
H_p^{a-1/p'+0}(\comega)$ for some $p$, cf.\ also {\rm (7.8)}.

Let $p>n/a$. Then
{\rm (7.5)} is
solvable when $f$ is in a subspace of $L_p(\Omega )$ with finite codimension, and the solutions satisfy
{\rm (7.9)ff}.

A similar statement holds for $f\in L_\infty (\Omega )$ with solutions
satisfying {\rm (7.10)}, and for $f\in C^{t+0}(\comega)$ with solutions
satisfying {\rm (7.11)}. Moreover, {\rm (7.12)} holds with Fredholm solvability.
\endproclaim

Since $a>0$, we can also apply Theorem 7.1 $3^\circ$ with $M=1$.
Recall that $\gamma _{a-1,0}u$ is a constant times $\gamma
_0(d(x)^{1-a}u)$. According to the theorem,
the nonhomogeneous Dirichlet problem
$$
r^+P_au=f, \quad \gamma _{0}d(x)^{1-a}u=\varphi ,\tag7.13
$$
is, when $s>a-1/p'$, Fredholm solvable for $f\in \ol H_p^{s-2a}(\Omega )$, $\varphi \in
B_p^{s-a+1-1/p}(\partial\Omega )$, with solution 
$u\in H_p^{(a-1)(s)}(\comega )$. 

Since $s>(a-1)+1-1/p'$, and $a-1>-1$,
Theorem 5.4 and its corollary  apply to show that when $s>n/p$, 
$$
\aligned
H_p^{(a-1)(s)}(\Omega )&\subset e^+d(x)^{a-1}\ol
H_p^{s-a+1}(\Omega )+\dot H^{s-0}_p(\comega)\\
&\subset e^+d(x)^{a-1}C^{s-a+1-n/p-0}(\comega)+\dot
C^{s-n/p-0}(\comega),
\endaligned\tag7.14
$$
where $-0$ can be left out except at certain values of $t$. 
(The $\dot C$-term is needed when $a<1$.)
Here we find:
$$
\multline
f\in L_p(\Omega ),\; \varphi \in  C^{a+1-1/p+0}(\partial\Omega )\implies 
\\
u\in e^+d(x)^{a-1}C^{a+1-n/p-0 }(\comega )\cap C^{2a-n/p -0}(\Omega
)+\dot C^{2a-n/p-0}(\comega),
\endmultline\tag7.15
$$
when $p>n/(a+1)$;
the $-0$ can be left out when $a-n/p$, $2a-n/p$ and $a-1/p$ are not integer.
For $p\to\infty $ this gives, since $L_\infty (\Omega )\subset
L_p(\Omega )$ and $C^{a+1}(\partial\Omega )\subset
C^{a+1-1/p+0}(\partial\Omega )$ for all $p$,
$$
f\in L_\infty (\Omega ),\; \varphi \in  C^{a+1}(\partial\Omega )\implies
u\in e^+ d(x)^{a-1}C^{a+1-0 }(\comega )\cap C^{2a-0}(\Omega )+\dot C^{2a-0}(\comega).
\tag7.16
$$
For $t\ge 0$ we likewise find
$$
\multline
f\in C^{t+0}(\comega),\; \varphi \in  C^{t+a+1}(\partial\Omega )\implies
\\
u\in e^+d(x)^{a-1}C^{t+a+1-0 }(\comega )\cap C^{t+2a-0}(\Omega )+\dot C^{t+2a-0}(\comega).
\endmultline\tag7.17
$$
In each of these situations, there is solvability when the data $\{f,\varphi \}$ 
are subject to a finite dimensional linear condition.
We recall moreover from Theorem 7.1 $3^\circ$ that 
$$
f\in C^\infty (\comega),\; \varphi \in  C^{\infty }(\partial\Omega )\iff
u\in e^+d(x)^{a-1}C^{\infty  }(\comega )\,\bigl(=\E_{a-1}(\comega)\bigr),\tag7.18
$$
with Fredholm solvability, when $u\in H_p^{(a-1)(s)}(\comega )$ for
some $s,p$ with $s>a-1/p'$.

We have then obtained:

\proclaim{Theorem 7.4} Hypotheses as in Theorem {\rm 7.3}. Consider
the nonhomogeneous Dirichlet problem {\rm (7.13)}.

Let $p>n/(a+1)$. For $u\in H^{(a-1)(\sigma )}(\comega)$ with $\sigma
>\max\{a-1/p',n/p\}$, cf.\ also {\rm (7.14)},  {\rm (7.13)} is
solvable when $f\in L_p(\Omega )$, $\varphi \in
C^{a+1-1/p+0}(\partial\Omega )$, subject to a finite dimensional
linear condition, with solutions satisfying
{\rm (7.15)ff}. 

A similar statement holds when $f\in L_\infty (\Omega )$, $\varphi \in
C^{a+1}(\partial\Omega )$, with solutions satisfying {\rm (7.16)},
and when $f\in C^{t+0}(\comega )$, $\varphi \in
C^{t+a+1}(\partial\Omega )$, with solutions satisfying {\rm (7.17)}.

Moreover, {\rm (7.18)} holds with Fredholm solvability. 
\endproclaim

Note that since $a$ can be any positive number, this covers powers between
0 and 1 of $\Delta ^2$, $\Delta ^3$, etc. When $a>1$, we can also
apply Theorem 7.1 $3^\circ$ for larger $M$ (namely for
$M<a+1$), which gives  natural extensions of Theorem 7.4. Details are
left to the reader. 

The theory moreover applies to $a$'th powers of
$2m$-order strongly elliptic differential operators, since they are of
order $2am$ and type $am$, and have
factorization index $am$, cf.\ Example 3.2.
The power $a$ can also be taken complex.

Other boundary operators (e.g.\ the Neumann operator $\gamma _{a-1,1}$
in lieu of $\gamma _{a-1,0}$ in (7.13), and more generally combinations of $\varrho _{\mu ,M}$ with
suitable $\psi $do's) can also be investigated, and one can make
applications to mixed problems and transmission problems, and to spectral
asymptotics. The solvability properties in H\"older spaces can be
sharpened slightly by applying the $\psi $do techniques
directly to scales of H\"older-Zygmund spaces $B^s_{\infty ,\infty }$.
We shall
return to these subjects in subsequent works.

\example{Remark 7.5}
The notes \cite{H65}, labeled Chapter II,  were given to me by Lars H\"ormander in
1980, but I have only studied them in depth recently. They have been
given to a number of people, but those colleagues that
I have asked (in order to find the missing Chapter I) have lost track of
them. I have typed  the text in \TeX{} (with comments on
misprints etc.),  and am willing to send it
to interested readers; it can also be found on my homepage.
\endexample


\Refs
\widestnumber\key{[GH90]}


\ref \no[AM09]\by P. Albin and R. B. Melrose\paper Fredholm realizations
of elliptic symbols on manifolds with boundary \jour J. Reine
Angew. Math. \vol 627 \yr 2009 \pages 155--181 \endref

\ref\no[B69]\by 
  L.~Boutet de Monvel  \paper Op\'erateurs  pseudo-diff\'erentiels
elliptiques et probl\graveaccent emes aux limites\jour Ann. Inst. Fourier
  \vol 19\yr 1969  
\pages  169--268 \endref

\ref\no[B71]\by 
  L.~Boutet de Monvel  \paper Boundary problems for pseudo-differential
operators\jour  
 {Acta Math.} \vol126\pages  11--51 \yr 1971\endref

\ref\no[B79] \by L. Boutet de Monvel\paper Lacunas and transmissions
\inbook Annals of Math. Studies \publaddr Princeton \vol 91 \yr 1979 \pages 209--218 \endref

\ref\no[BGR10]\by K. Bogdan, T. Grzywny and M. Ryznar\paper Heat
kernel estimates for the fractional Laplacian with Dirichlet
conditions \jour Ann. of Prob.\yr 2010 \vol38 \pages 1901--1023
\endref


\ref\no[CD01] \by O. Chkadua and R. Duduchava\paper Pseudodifferential
equations on manifolds with boundary: Fredholm property and
asymptotics
\jour Math. Nachr.\yr2001 \vol 222 \pages 79--139
\endref

\ref\no[E81]\by G. Eskin\book Boundary value problems for elliptic
pseudodifferential equations \publ Amer. Math. Soc. \publaddr
Providence, R.I.\yr 1981
 \endref

\ref\no[FG14]\by R. Frank and L. Geisinger\paper Refined semiclassical
asymptotics for fractional powers of the Laplace operator\finalinfo
arXiv:1105.5181, to appear in J. Reine Angew. Math
\endref

\ref\no[F86]\by J. Franke \paper Elliptische Randwertprobleme in
Besov-Triebel-Lizorkin-Ra\"umen \finalinfo Dissertation,
Friedrich-Schiller-Universit\"at Jena \yr1986 \endref

\ref\no[GMS09]\by M. Gonzalez, Rafe Mazzeo and Y. Sire \paper Singular
solutions of fractional order conformal Laplacians \jour
J. Geom. Anal. \vol 22 \yr 2012\pages 845-–863
\endref

\ref\no[G90] \by G. Grubb \paper Pseudo-differential boundary problems
in $L_p$-spaces \jour Comm. Part. Diff. Eq. \vol 13 \yr 1990 \pages
289--340
\endref




 \ref\no[G96]\by 
{G.~Grubb}\book Functional calculus of pseudodifferential
     boundary problems.
 Pro\-gress in Math.\ vol.\ 65, Second Edition \publ  Birkh\"auser
\publaddr  Boston \yr 1996\finalinfo first edition issued 1986\endref

\ref\no[G09]\by G. Grubb\book Distributions and operators. Graduate
Texts in Mathematics, 252 \publ Springer \publaddr New York\yr 2009
 \endref


\ref\no[GH90] \by G. Grubb and L. H\"ormander \paper The transmission
property
\jour Math. Scand.  \vol 67 \yr 1990 \pages
273--289
\endref

\ref\no[GK93]\by G. Grubb and N. J. Kokholm \paper A global calculus
of parameter-dependent pseudodifferential boundary problems in $L_p$
Sobolev spaces \jour Acta Math. \yr 1993 \vol 171 \pages 165--229 \endref



\ref\no[HS08]
\by G. Harutyunyan and B.-W. Schulze \book Elliptic mixed,
transmission and singular crack problems. EMS Tracts in Mathematics, 4
\publ European Mathematical Society (EMS) \publaddr Z\"urich \yr  2008
\endref

\ref\no[HP79]\by A. Hirschowitz and A. Piriou\paper  Propri\'e{}t\'e{}s
de transmission pour les distributions int\'e{}grales de Fourier \jour Comm.
Part. Diff. Eq. \vol 4 \yr1979\pages 113-–217
\endref

\ref\no[H63]\by L. H\"o{}rmander\book Linear partial differential operators\yr1963 
\endref

\ref\no[H65]\by L. H\"o{}rmander\book Ch.\ II, Boundary problems for
``classical'' pseudo-differential operators \finalinfo photocopied lecture notes
at Inst. Adv. Study, Princeton\yr1965 
\endref



\ref\no[H83]\by L. H\"o{}rmander\book The analysis of linear partial
differential operators, I \publ Springer Verlag\publaddr Berlin, New
York\yr 1983 
 \endref

\ref\no[H85]\by L. H\"o{}rmander\book The analysis of linear partial
differential operators, III \publ Springer Verlag\publaddr Berlin, New
York\yr 1985
 \endref

\ref \no[M93]\by R. B. Melrose\book The Atiyah-Patodi-Singer index
theorem \publ A. K. Peters \publaddr Wellesley, MA \yr 1993 \endref


\ref\no[RS84]\by S. Rempel and B.-W. Schulze \paper Complex powers for
pseudo-differential boundary problems II \jour Math. Nachr. \vol 116
\yr 1984 \pages 269--314\endref

\ref\no[R96]\by Y. Roitberg \book
Elliptic boundary value problems in the spaces of distributions.
Mathematics and its Applications \vol 384 \publ Kluwer Academic
Publishers Group\publaddr Dordrecht \yr 1996 \pages 415 pp\endref

\ref\no[RS69]\by Y.A. Roitberg and V. Sheftel \paper
A homeomorphism theorem for elliptic systems, and its applications
\jour Mat. Sb. (N.S.) \vol 78 \yr 1969 \pages 446-–472
\endref



\ref\no[RS14] \by X. Ros-Oton and J. Serra \paper The Dirichlet
problem for the fractional Laplacian
\jour  J. Math. Pures Appl.  
\yr 2014 \pages 275-302 \vol 101
\endref



\ref\key[S67]\by R. T. Seeley \paper Complex powers of an elliptic
operator\jour Amer. Math. Soc. Proceedings Symposia Pure Math. \vol
10\yr 1967\pages 288--307
\endref


\ref\no[S94]\by E. Shargorodsky\paper An $L_p$-analogue of the
Vishik-Eskin theory \inbook Memoirs on Differential Equations and
Mathematical Physics, Vol. 2\publ Math. Inst.
Georgian Acad. Sci. \publaddr Tblisi\yr 1994\pages 41--146
\endref

\ref\no[T95]\by H. Triebel \book Interpolation theory, function
spaces, Differential operators (2nd edition)\publ J. A. Barth
\publaddr Leipzig \yr 1995
\endref

\ref\no[VE65] \by M. I. Vishik and G. I. Eskin \paper
Convolution equations in a bounded region
\jour Uspehi Mat. Nauk \vol 20 \yr 1965 \pages 89--152\transl\nofrills
English translation in \jour Russian Math. Surveys \vol 20 \yr1965
\pages 86--151\endref


\ref\no[VE67] \by M. I. Vishik and G. I. Eskin \paper
Convolution equations of variable order
\jour Trudy Mosk. Mat. Obsc. \vol 16 \yr 1967 \pages 25--50
\transl\nofrills  English translation in \jour Trans. Moscow
Math. Soc. \vol 16 \yr 1967 \pages 27--52\endref

\endRefs
\enddocument